\documentclass[3p]{elsarticle}
\usepackage[utf8]{inputenc}
\usepackage{geometry}
\usepackage{amsmath}
\usepackage{cancel}
\usepackage{bbm}	
\usepackage{graphicx}
\usepackage{placeins}
\usepackage{color,soul, colortbl}
\usepackage{mathtools}
\usepackage{commath}
\usepackage{empheq}
\usepackage{tabularx}
\usepackage{cases}
\usepackage{amssymb}	
\usepackage{siunitx} 
\usepackage{tikz, tikz-3dplot}

\date{}

\def\bfpsi{\boldsymbol{\psi}} 
\def\bfR{{\bf R}}  
\def\bfV{{\bf V}}  
\def\bfBXV{\hat{\bf B}_{V  x}}  
\def\bfBYV{\hat{\bf B}_{V  y}}  
\def\bfBXE{\hat{\bf B}_{E  x}}  
\def\bfBYE{\hat{\bf B}_{E  y}}  
\def\bfE{{\bf E}}  
\def\bEE{\boldsymbol{\mathcal{E}}}  
\def\EE{\mathcal{E}}
\def\bfI{{\bf I}}  
\def\bfD{{\bf D}}  
\def\bfu{{\bf u}}  
\def\bfh{{\bf h}}  
\def\bfn{{\bf n}}  
\def\bfne{{\bfn}_e}  
\def\ww{{ w}}  
\def\bfw{{\bf w}}  
\def\bfwu{{{\bf w}_V  }}  
\def\bfwe{{{\bf w}_E }}  
\def\bfe{{\bf e}}  
\def\bfx{{\bf x}}  
\def\bfI{{\bf I}}  
\def\bfA{{\bf A}}  
\def\bfQ{{\bf Q}}  
\def\bfzero{{\bf 0}} 
\def\div{{\text{div}\,}} 
\def\dd{{\text{d}}}  
\def\bfr{{\bf r}}  
\def\htt{{t}}    
\def\htx{{x}}    
\def\htH{{H}}    
\def\hth{{h}}    
\def\htrho{{\rho}} 
\def\bfU{{\bf U}}  
\def\bfUi{{\bfU}^{-} } 
\def\bfEi{{\bfE}^{-} } 
\def\bfwi{{\bfw}^{\,\,\,-}} 
\def\bfUe{{\bfU}^{+} } 
\def\bfEe{{\bfE}^{+} } 
\def\bfVi{{\bfV}^{-} } 
\def\bfVe{{\bfV}^{+} } 
\def\bfF{{\bf F}}  
\def\bfG{{\bf G}}  
\def\bfQi{\bfQ^{-}}  
\def\bfQe{\bfQ^{+}}  
\def\nbfF{\overline{\bfF}}  
\def\nbfA{\overline{\bfA}}  
\def\bfGi{\bfG^{-}}  
\def\bfGe{\bfG^{+}}  
\def\nbfG{\overline{\bfG}}  
\def\nbfQ{\overline{\bfQ}}  
\def\Heav{\mathcal{H}}

\def\vel{\bar{u}}
\def\bvel{\bar{\bf u} }
\def\vel{\bar{u} }

\def\bfnabla{\boldsymbol{\nabla}} 
\def\partialt{\partial_{\htt}} 
\def\bfS{{\bf S}}  
\def\nel{n_{el}}  
\def\trace{\text{tr}} 
\def\sigmao{\sigma_o} 

\def\bfsigma{\bar{\boldsymbol{\sigma}}} 
\def\bfsigmaOne{\bar{\boldsymbol{\sigma}}_{N}} 
\def\bfsigmaTwo{\bar{\boldsymbol{\sigma}}_{B}} 

\def\veli{{\vel}^{-} } 
\def\vele{{\vel}^{+} } 
\def\hi{{\hth}^{-} } 
\def\he{{\hth}^{+} } 

\def\bfsigmaTwoz{\bar{\boldsymbol{\sigma}}_{B }} 
\def\bfsigmaTwoa{\bar{\boldsymbol{\sigma}}_{B(1) }} 
\def\bfsigmaTwob{\bar{\boldsymbol{\sigma}}_{B(2) }} 
\def\bfsigmaTwoc{\bar{\boldsymbol{\sigma}}_{B(3) }} 

\def\Lh{$L^\infty_h$}
\def\Lvel{$L^\infty_{\vel}$}
\def\Lth{$L^2_h$}
\def\Ltvel{$L^2_{\vel}$}

\newcommand{\da}{\mbox{d}\Omega}
\newcommand{\dl}{\mbox{d}\Gamma}
\newcommand{\Real}{\mathbb{R}}     

\newcommand{\normc}[1]{\norm{#1}_c \,} 
\newcommand{\vect}[1]{\text{vec}\left({#1} \right)} 
\newcommand{\hbfpsi}[1]{\hat{\bfpsi}_{#1} } 

\newcommand{\signf}[1]{\text{sign}{ \left(#1 \right)} \,} 

\newcommand{\derive}[2]{ \frac{\partial #1 }{\partial #2}}
\newcommand{\derivev}[1]{ \derive{#1}{\bfV_e} }
\newcommand{\derivee}[1]{ \derive{#1}{\bfE_e} }
\newcommand{\derivees}[1]{ \derive{#1}{\EE} }
\def\bfNV{\hat{\bf N}_V}  
\def\bfNE{\hat{\bf N}_E}  

\definecolor{Gray}{gray}{0.9}

\begin{document}

	\begin{frontmatter}
		
		\title{ Regularized Approach for Bingham Viscoplastic Shallow Flow Using the Discontinuous Galerkin Method }
		
		\author[mymainaddress]{{Felipe Fern\'{a}ndez}\corref{mycorrespondingauthor}}
		\cortext[mycorrespondingauthor]{Corresponding author}
		\ead{felipe.fernandez.ayala@gmail.com}
		
		\author[mymainaddress]{Sof\'{i}a L\'{o}pez-Ord\'{o}ñez}
            \ead{sofia.lopezo@epn.edu.ec}
		
		\author[mymainaddress]{Sergio Gonz\'{a}lez-Andrade}
            \ead{sergio.gonzalez@epn.edu.ec}

		\address[mymainaddress]{Research Center on Mathematical Modeling (MODEMAT), Escuela Politécnica Nacional, Quito, Ecuador.}
		
		\begin{abstract}
                This paper aims to simulate viscoplastic flow in a shallow-water regime.
                We specifically use the Bingham model in which the material behaves as a solid if the stress is below a certain threshold, otherwise, it moves as a fluid.
                The main difficulty of this problem is the coupling of the shallow-water equations with the viscoplastic constitutive laws and the high computational effort needed in its solution.
                Although there have been many studies of this problem, most of these works use explicit methods with simplified empirical models.
                
                In our work, to accommodate non-uniform grids and complicated geometries, we use the discontinuous Galerkin method to solve shallow viscoplastic flows.
                This method is attractive due to its high parallelization, h- and p-adaptivity, and ability to capture shocks.
                Additionally,  we treat the discontinuities in the interfaces between elements with numerical fluxes that ensure a stable solution of the nonlinear hyperbolic equations. 
                To couple the Bingham model with the shallow-water equations, we regularize the problem with three alternatives.
                Finally, in order to show the effectiveness of our approach, we perform numerical examples for the usual benchmarks of the shallow-water equations.
		\end{abstract}
		
		
	\end{frontmatter}

\section{Introduction}

The shallow-water equations accurately model flows in channels, rivers, seas, floods, and tsunamis \cite{Costa2005, Wesseling2009}.
These equations are derived from the Navier-Stokes equations assuming that the fluid height is much smaller than the horizontal length scale \cite{Costa2005}.
In this procedure, the continuity and momentum equations are depth-integrated.
The resulting equations present uniform horizontal velocities, a hydrostatic pressure distribution, and neglected vertical velocities.
As such, the shallow-water equations govern flows with long waves or shallow fluids.
Since we are interested in modeling geophysical flows where the horizontal length scale (kilometers) is much larger than the vertical one (meters), the shallow-water equations provide us with a suitable framework.
The shallow-water equations were first introduced by De Saint Venant in 1864 and Boussinesq in 1872 \cite{Costa2005}, and researchers have solved them to model geophysical flows producing many methods and countless papers on these topics.
We do not attempt to give an exhaustive survey of the literature in this paper.
However, we highlight some relevant references to provide a discussion of the state-of-the-art, and we describe our contribution made with this work.

Researchers have modeled lava and avalanches using the traditional finite difference method and explicit upwind schemes.
In the pioneering work by Savage and Hutter \cite{Savage1989}, the authors solved depth-averaged equations similar to the shallow-water ones and use a Mohr-Coulomb constitutive relation to model the motion of granular material through a rough incline in one dimension. 
They use finite differences and reported difficulties with the Eulerian approach, and claim that the Lagrangian approach is simpler, efficient, and reliable.
Additionally, the authors validated the numerically-predicted motion of the granular flow  with laboratory experiments.
This work was extended to two dimensions in \cite{Hutter1993}, to pore fluid in \cite{Iverson1997}, and to complex topography in \cite{Gray1999}.
Unfortunately, classical finite difference methods can be expected to break down near discontinuities in the solution where the differential equation does not hold \cite{Leveque2002}.
Shock-capturing methods for finite differences have been proposed to resolve this issue.
For instance, in \cite{Tai2002}, researchers use an Eulerian shock-capturing non-oscillatory central scheme which is built upon the classical Lax-Friedrichs scheme and combined with a front-tracking method.
Other authors \cite{Kelfoun2005,  Kelfoun2009, Kelfoun2016} used a double upwind Eulerian scheme to capture shocks.
Besides the shock-capturing difficulties, the finite difference method requires extending the computational stencil to achieve higher-order accuracy, and this might be problematic for two-dimensional problems and complex geometries.
Furthermore, accurate and stable treatment of this extended stencil at the boundaries is not trivial.

Similarly, the finite volume method has been used to solve the shallow-water equations \cite{Leveque2002}.
Applications of the finite volume method used to model geophysical flows like lava, tsunamis, granular flow, and avalanches have been successfully treated in \cite{Denlinger2001,Pitman2003, Costa2005b, Patra2005, Gallardo2007}.
Even though essentially non-oscillatory schemes used in the finite volume method allow the use of nonuniform or unstructured grids, their computational cost and complexity are high \cite{Hesthaven2017}.
Also, an important difficulty of the numerical solutions is the proper tracking of wet/dry fronts.
Most of the finite volume methods in the literature use explicit schemes, so the Courant-Friedrichs-Lewy (CFL) condition controls the stability.
The CFL condition is determined with the maximum characteristic speed of the governing equations, which depends on the square root of the flow thickness. 
Consequently, numerical schemes that produce negative thickness would experience numerical issues and  are additionally physically unrealistic.
To solve this problem, a simple but effective method is to set any negative thickness to zero after each time iteration \cite{Bi2014}. 
However,  this method might not preserve the total mass of the system. 

Real applications of the shallow-water equations, e.g., geophysical flow, are often dominated by source terms such as bed topography, friction, etc.
The presence of source terms in hyperbolic conservation laws often admits steady-state solutions where the fluxes and source terms must balance each other. 
Researchers introduced the so-called well-balanced methods to satisfy this balance.
In the case of the shallow-water equations, a well-balanced method preserves still water at rest in a steady state \cite{Bermudez1994}.
A straightforward approach may fail to preserve exactly this steady state and introduce spurious oscillations near the steady state.
Another advantage of well-balanced methods is that they can solve small perturbations of such steady-state solutions with coarse meshes.
To preserve stationary steady-state solutions, in \cite{Kurganov2018}, the authors adopted a central-upwind finite-volume scheme which also preserves the positivity of the height.
Other conservative methods to treat dry/wet fronts have been proposed, e.g., \cite{Begnudelli2007, Liang2009}.

Besides the computational instabilities described above due to numerical modeling of wet/dry fronts, other computational issues might appear when friction or material modeling terms depend on the flow thickness or its inverse, e.g., Manning formulation \cite{Bi2014}.
For instance, in \cite{Bi2014,Begnudelli2007}, the friction terms are solved by a semi-implicit scheme to prevent such instabilities. 
Similarly, implicit-explicit Runge-Kutta methods have been proposed \cite{Russo2005} to deal with stiff source terms.
Another advantage of using semi-implicit or implicit-explicit schemes is that the time step is not restricted by the CFL condition.

Besides the finite differences and finite volume method, the finite element method has the potential for high-order accuracy and ensures geometry flexibility. 
However, this method is inherently implicit since it requires the inverse of the mass matrix.
Consequently, the finite element method for time-dependent problems can be computationally more expensive than explicit finite differences and volume methods.
Additionally, the finite element method cannot exploit upwinding in a straightforward manner since the trial/basis functions depend solely on the grid \cite{Hesthaven2017}.
Also, straightforward use of equal order approximating spaces for height and velocity in the equations can lead to spurious spatial oscillations \cite{Lynch1979}.

Another option is using the discontinuous Galerkin method, which combines the advantages of both: the finite element and finite volume methods.
It accommodates nonuniform grids and complex geometries and takes advantage of higher-order methods.
As the finite element method, enriching the local basis enables higher-order accuracy to the discontinuous Galerkin method.
And similarly to the finite volume method, the discontinuous Galerkin method uses numerical fluxes at element interfaces, allowing the exchange of information.
The computation of such fluxes is vital for the accuracy and stability of the numerical solution.
This method is attractive due to its high parallelization, h- and p-adaptivity, and ability to capture shocks \cite{Hesthaven2017}.
These benefits are not free, and there is an additional computational cost since decoupling elements increases the total number of degrees of freedom.
However, the mass matrix is local rather than global and can be inverted at a low cost, so explicit schemes are computationally inexpensive.
Also, higher-order polynomial representations lead to artificial oscillations, which tend to be local and at discontinuities of the solution. 
Even though these oscillations do not destroy the accuracy of the method, these oscillations are spurious, and in countless situations, they cannot be tolerated.
Fortunately, slope-limiting techniques and other approaches \cite{Cockburn1989} have been proposed to reduce the effect of these undesired oscillations.
For these reasons, the discontinuous Galerkin method is a usual choice to solve hyperbolic equations.

The discontinuous Galerkin method was first proposed by \cite{Reed1973} to solve the steady-state neutron transport equation.
Extensions to solve hyperbolic systems were treated in \cite{Brezzi2004, Ern2006}.
Numerous applications of this method have been developed, and many textbooks offer different introductions to this active research area \cite{Hesthaven2017}. 
Researchers have successfully used the discontinuous Galerkin method to solve the shallow-water equations, e.g., \cite{Schwanenberg2000, Li2001, Aizinger2002, Giraldo2002, Eskilsson2004, Dawson2005, Dawson2006, Navas2020}.
Also, many studies included the treatment of wet/dry fronts, e.g., \cite{Bokhove2005, Ern2008, Bunya2009, Xing2010, Khan2014, Gerhard2015, Caviedes2020}.
Usually, these studies use explicit time stepping and are restricted by the CFL condition.

The rheological behavior of geophysical materials like lava, snow avalanches,  and debris flows is difficult to determine because the material is heterogeneous and includes irregular particles over a wide range of sizes \cite{Ancey2007}. 
Additionally, researchers often use data from past events over complex mountain topographies or run small-scale experiments to validate these rheological models.
Therefore, this validation is not straightforward.
For these reasons, the understanding of the behavior of geophysical flows is in active debate and continuous research.
Despite these difficulties, geophysical flows are often modeled as plastic materials, i.e., the material yields and start to flow once its stress state is above a certain critical threshold.
Two plastic theories are commonly used in fluid dynamics: Coulomb plasticity and viscoplasticity \cite{Ancey2007}.
A review paper about the constitutive rheological behavior of geophysical flows can be found in \cite{Ancey2007}, and for a review about numerical simulations of viscoplastic flows, the reader is referred to \cite{Saramito2017}.

In this paper, we are interested in modeling geophysical fluids as viscoplastic materials and specifically, we use the Bingham model.
In this model, if the material's stress magnitude is below a certain critical threshold, i.e., the material's yield stress, there is no deformation, and the material behaves like a solid.
However, if the yield stress is exceeded there is deformation, and the material flows like a fluid.
So, the materials' yield stress would be the stress at which the solid-state first starts to deform continuously, i.e., flows.
In the literature, we find studies of this problem that introduce the Bingham model in the shallow-water equations with simplified empirical models.
For instance, Coussot \cite{Coussot1994} developed a reduced model for viscoplastic materials and validated it on steady uniform flows of muddy debris flows in a laboratory flume. 
This friction expression was introduced in a numerical model that uses the finite-volume method to solve the shallow-water equations \cite{Laigle1997}.
In \cite{Kelfoun2016}, researchers use a model function of the velocities (instead of the strain-rate tensor). 
The numerical solution of this model is somewhat easy, but the model does not accurately capture the nature of the Bingham material.

To rigorously include the Bingham model in the shallow-water equations, researchers \cite{Bresch2009, Ionescu2010} use depth-integration of the three-dimensional equations and obtain the constitutive relation between in-plane depth-averaged stresses with in-plane rate deformations. 
It is important to note that the resulting constitutive law for the Bingham model for the shallow-water equations has a similar structure to the three-dimensional one.
The main difficulty of this problem is the coupling of the shallow-water equations with the viscoplastic constitutive laws and the high computational effort needed in its solution.
The plastic stress term of this model is nondifferentiable, and the variational form of the shallow-water equations coupled with the Bingham model results in an inequality for the momentum equation, which is difficult to solve efficiently.
To circumvent this difficulty, researchers use regularization or the augmented Lagrangian method.
Ionescu in \cite{Ionescu2010} used the implicit (backward) Euler scheme and solved the variational inequality of the momentum equation with iterative decomposition-coordination formulation coupled with the augmented Lagrangian method using the finite element discretization. 
The author uses the finite volume method with an upwind strategy to solve the continuum equation.
Additionally, the author deduced and included frictional contact of the viscoplastic material with the plane slope.
This work, developed for plane slopes, was extended to general topographies in \cite{Ionescu2013a, Ionescu2013b}.
However, it does not consider well-balanced properties or wet/dry front treatment.
In \cite{Bresch2009}, researchers use the augmented Lagrangian method with well-balanced properties coupled with the finite volume discretization to solve one-dimensional Bingham viscoplastic flow.
In \cite{Acary2012}, researchers extended the previous work in \cite{Bresch2009} and derived an integrated Herschel-Bulkley model (which generalized the Bingham law) for the shallow-water equations and solved this problem.
In an additional extension, another duality method namely Berm\'udez-Moreno was implemented to model viscoplastic avalanches \cite{Fernandez2014} instead of the augmented Lagrangian method.
In this work, the authors also include the treatment of dry/wet front.
They reported that the Berm\'udez-Moreno approach is in general more efficient than the augmented Lagrangian one.
This work was extended to two-dimensional problems in \cite{Fernandez2018}.
Shallow viscoplastic models coupled with thermal cooling were
studied in \cite{Balmforth2004} for lava domes and in \cite{Bernabeu2016} for lava flows on complex tridimensional topographies.


Instead of using the augmented Lagrangian or the Berm\'udez-Moreno methods to solve shallow viscoplastic models, we propose to use a simple and easy-to-implement regularization approach as in \cite{DeLosReyes2012a}.
We present and study three alternatives of regularization, and demonstrate their great performance in numerical tests. 
Two of these regularizations are of the local type and one is global.
The local type regularizations enforce a local smoothing where the stress results in a $C^1$ function, and in contrast the global type regularization is $C^\infty$ smooth. 
The idea is that we modify the problem and recover standard equations for a nonlinear Newtonian fluid.
With the help of these regularizations, the solid region is modeled as fluid with large viscosity, and the numerical solution of the problem allows us to use well-known methods such as Newton-Raphson to obtain quadratic convergence.
Also, to accommodate non-uniform grids and complicated geometries, we use the discontinuous Galerkin method.
As mentioned, this method is attractive due to its high parallelization, h- and p-adaptivity, and ability to capture shocks.
We treat the discontinuities in the interfaces between elements with numerical fluxes that ensure a stable solution to the problem. 
In Sections \ref{Sec:ShallowWaterEquations} and \ref{Sec:VectorialEquations}, we present the shallow-water equations and their vectorial form.
Section \ref{Sec:Regularization} details our regularization approach, and we describe the discontinuous Galerkin method in Section \ref{Sec:DG}.
Numerical examples are performed in Section \ref{Sec:Examples} for the usual benchmarks showing the effectiveness of our approach, and concluding remarks are drawn in Section \ref{Sec:Conclusions}.

\section{The shallow-water equations}
\label{Sec:ShallowWaterEquations}
\begin{figure}[!ht]
\centering
	\begin{tikzpicture} [scale=1]
    	\fontsize{9pt}{9pt}\selectfont
    	\def\xo{-.25cm}		
    	\def\yo{-.25cm}		
    	\draw[draw=none, use as bounding box](0cm+\xo,0cm+\yo) rectangle (8cm+\xo,5cm+\yo);
    	
    	\begin{scope} [rotate=20]
        	\def\xyo{0}
        	\def\dxyo{-.2}
        	\def\dxa{7.5}	
    	    \def\dya{3.5}
    	    \draw[->] (\xyo,\dxyo)--(\xyo,\xyo+\dya);
    	    \draw[->] (\dxyo,\xyo)--(\xyo+\dxa,\xyo);
        	\node at (\xyo+\dxa,\xyo-.3) {$x$};
        	
        	\draw[ultra thick, orange] plot [smooth,tension=0.7] coordinates {(0,2.25) (1.5,1) (3,1.5) (4.5,0.5) (6,1) (7,0)};
        	\node[orange] at (\xyo+\dxa-.5,\xyo+.9) {$H(t,\bfx)$};
        	
        	\draw[ultra thick, magenta] plot [smooth,tension=0.7] coordinates {(0,3) (1.5,2.5) (3,2.75) (4.5,2.5) (6,2.75) (7,2.25)};
        	\node[magenta] at (\xyo+\dxa-.5,\xyo+.75+2) {$\xi(t,\bfx)$};
        	
        	\draw[thick,cyan,<->] (6,1)--(6,2.75);
        	\node[cyan] at (5.5,1.6) {$h(t,\bfx)$};
        	
        	\def\ud{.75}	
    		\draw[thick,red,-latex] (\xyo,\xyo)--(\xyo+\ud,\xyo) node[above] {$\bfe_1$};
    		\draw[thick,red,-latex] (\xyo,\xyo)--(\xyo,\xyo+\ud) node[left] {$\bfe_3$};
    		\node at (\dxyo,\dxyo) {$\bf 0$};
		\end{scope}
		\draw[dashed] (0,0)--(5,0);
		\draw[thin,gray] (2,0) arc (0:20:2);
		\node at (1.5,.25) {$\alpha$};
		
    	\def\dxyo{-.2}
    	\def\dxa{7.5}	
	    \def\dya{3.5}
	    \def\ud{.75}	
		\draw[thick,blue,-latex] (0,0)--(\ud,0) node[below] {$\bfe_1'$};
		\draw[thick,blue,-latex] (0,0)--(0,\ud) node[above] {$\bfe_3'$};
		
    \end{tikzpicture}
    
    \caption{Fluid and terrain sketch for the shallow-water equations.}
	\label{Fig:shallow}
\end{figure}
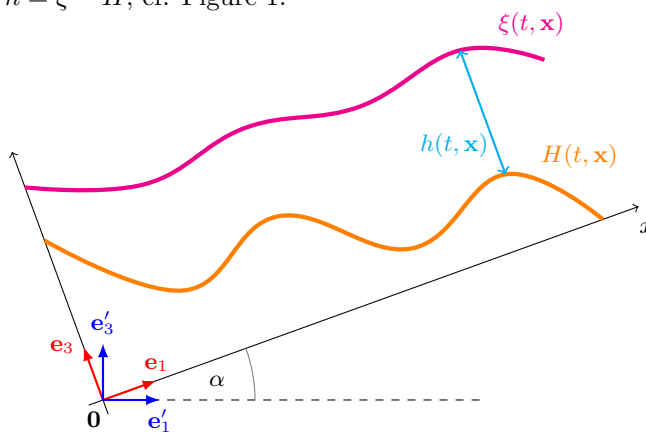

Consider a three-dimensional domain where a fluid with a free surface is under the influence of gravity acting in the opposite direction of $\bfe_3'$.
The fluid elevation is $\xi(t, x_1, x_2)$ and the terrain surface is $H(t, x_1,x_2)$, so the fluid depth is $h=\xi - H$, cf. Figure \ref{Fig:shallow}.

The shallow-water assumption allows us to approximate the horizontal velocity components $u_1$ and $u_2$ by its depth-average $\vel_1$ and $\vel_2$ \cite{Wesseling2009}, i.e.,
\begin{equation}
    \vel_\alpha(t, x_1,x_2) = \frac{1}{h(t,x_1,x_2)} \int_H^\xi u_\alpha(t, \bfx) \, \dd x_3 \,, \quad \alpha = 1, 2\,.
    \label{Eq:averageU}
\end{equation}

Assuming the density $\rho$ is constant, we integrate the continuity equation $\div \bfu = 0$ over $x_3$, replace the boundary conditions, and use the Leibniz integral rule, where the resulting first shallow-water equation, i.e., the continuity equation, is
\begin{eqnarray}
    0
    &=& \frac{\partial h}{\partial t} + \frac{\partial}{\partial x_1} \left( h \, \vel_1\right) + \frac{\partial}{\partial x_2} \left( h \, \vel_2\right) \,, 
\end{eqnarray}
where $t\in \Real^+$ is the time variable.

Similarly, the depth-averaged momentum equations are written as
\begin{eqnarray}
     \frac{\partial (h\, \vel_1)}{\partial t} 
    + \frac{\partial (h\, \vel_1^2)}{\partial x_1}
    +  \frac{\partial (h\,\vel_2\,\vel_1)}{\partial x_2} 
    + \frac{1}{2} g \, \cos(\alpha) \frac{\partial (h^2 ) }{\partial x_1}
    + g \,h\, \sin(\alpha) + g \, h \, \cos(\alpha) \frac{\partial H }{\partial x_1} 
    &=& \frac{T_1}{\rho}
    \\
     \frac{\partial (h\, \vel_2)}{\partial t} 
    +  \frac{\partial (h\,\vel_2\,\vel_1)}{\partial x_1}
    + \frac{\partial (h\, \vel_2^2)}{\partial x_2}
    + \frac{1}{2} g \, \cos(\alpha) \frac{\partial (h^2 ) }{\partial x_2}
    + g \,h\, \cos(\alpha) \frac{\partial H }{\partial x_2}
    &=& \frac{T_2}{\rho} \,,
\end{eqnarray}
where $g$ is the gravity and 
\begin{eqnarray}
    T_\alpha 
    &=&
    \frac{\partial \left( h \,  \bar{\sigma}_{\alpha1}  \right)}{\partial x_1}  
    +
    \frac{\partial \left( h \,  \bar{\sigma}_{\alpha2}  \right)}{\partial x_2} 
    +
    \tau_{s\alpha} - \tau_{b\alpha}\,,
\end{eqnarray}
where $\bar{\sigma}_{\alpha1}$ and $\bar{\sigma}_{\alpha2}$ are the depth-averaged stress components, and $\tau_{s\alpha}$ and $\tau_{b\alpha}$ are the shear stress tangent to the fluid surface and the bottom respectively.
In this study, we assume these surface stresses are zero, i.e, $\tau_{s\alpha}=0$ and $\tau_{b\alpha}=0$.

In the literature, we find several models for stress terms \cite{Costa2005,Dade1998,Kelfoun2009,Wesseling2009,Kelfoun2016}. 
For instance, the retarding stress to model Bingham flow in \cite{Kelfoun2016} is defined as
\begin{eqnarray}
    T_\alpha = \sigma_o\frac{\vel_\alpha}{|\vel|} + 3 \eta \frac{\vel_\alpha}{h}\,,
    \label{Eq:KelfounBingham}
\end{eqnarray}
where $\sigma_o$ is the yield stress, and $\eta$ the dynamic viscosity. 
As discussed in the introduction, in this model, the stress is a function of the velocity and not the strain-rate tensor, so the computational solution is less cumbersome.
However, these models are simplified empirical approximations and do not capture the true nature of the Bingham model.
For these reasons, we do not adopt them in our study.

In this work, we adopt a \emph{depth-integrated} Bingham law used to model viscoplastic avalanches in \cite{Bresch2009,Fernandez2010,Fernandez2014,Fernandez2018}. 
The depth-averaged stress in this formulation \cite{Fernandez2018} is the following 
\begin{eqnarray}
    \begin{cases}
       \bfsigma =2\eta\left( \bfD + \trace (\bfD)\,\bfI \right) 
       + \sqrt{2}\sigmao \, \frac{ \bfD + \trace (\bfD)\,\bfI}{\normc{\bfD }} & \text{if} \quad \normc{\bfD}\neq 0 \\
       \normc{\bfsigma} \leq \sqrt{2} \sigmao & \text{if} \quad \normc{\bfD}= 0\,,
    \end{cases}
    \label{Eq:bingham}
\end{eqnarray}
where $\bvel=[\,\vel_1\,,\,\vel_2]$,  $\eta$ is the viscosity, $\sigmao$ is the yield stress (plasticity threshold), $\bfI$ is the $2\times2$ identity, $\bfsigma$ is a symmetric second-order tensor of two-dimensions with components $\bar{\sigma}_{11}$, $\bar{\sigma}_{12}$,         $\bar{\sigma}_{21}$, $\bar{\sigma}_{22} $, 
and
\begin{eqnarray}
    \bfD = \frac{1}{2} \left( \nabla \bvel + \nabla \bvel^{\top} \right)\,,
    \label{Eq:Dtensor}
\end{eqnarray}
where $\left(\nabla \bvel \right)_{ij}=  \partial \vel_j/ \partial x_i$ for $i=1,2$ and $j=1,2$.
Lastly, the norm is defined as
\begin{eqnarray}
    \normc{p} = \sqrt{ \sum_{i=1}^2\sum_{j=1}^2 p_{ij}^2 + \left( \sum_{i=1}^2 p_{ii} \right)^2} \,.
\end{eqnarray}

\section{Vectorial equations}
\label{Sec:VectorialEquations}
Let $\Omega\in \Real^d $ be the fluid domain of space dimension $d = 1, 2$ with boundary $\Gamma$.
The shallow-water equations in the time interval  $t\in[0, t_f]$ are:
\begin{eqnarray}
    \frac{\partial \hth}{\partial \htt} + \frac{\partial (\hth \vel_1) }{\partial \htx_1} + \frac{\partial (\hth \vel_2)}{\partial \htx_2}  &=& 0 
    \label{Eq:shallow1}
    \\
    \frac{\partial (\hth \vel_1)}{\partial \htt} 
    + \frac{\partial (\hth \vel_1^2 )}{\partial \htx_1}
    + \frac{1}{2}g_c \frac{\partial (\hth^2 ) }{\partial \htx_1}
    +  \frac{\partial (\hth \vel_1\,\vel_2  )}{\partial \htx_2} 
    &=& 
    \frac{\partial \left( h \,  \bar{\sigma}_{11} /\htrho \right)}{\,\partial x_1}  
    +\frac{\partial \left( h \,  \bar{\sigma}_{12} / \htrho  \right)}{\partial x_2} 
    -g_s\,\hth - g_c\,\hth\, \frac{\partial \htH }{\partial \htx_1}
    \label{Eq:shallow2}
    \\
     \frac{\partial (\hth \vel_2)}{\partial \htt} 
    +  \frac{\partial (\hth \vel_1\,\vel_2 )}{\partial \htx_1}
    + \frac{\partial (\hth \vel_2^2)}{\partial \htx_2}
    + \frac{1}{2}g_c \frac{\partial (\hth^2 ) }{\partial \htx_2}
    &=& 
    \frac{\partial \left( h \,  \bar{\sigma}_{21} /\htrho \right)}{\,\partial x_1}  
    +
    \frac{\partial \left( h \,  \bar{\sigma}_{22} /\htrho \right)}{\,\partial x_2}
    - g_c\,\hth\, \frac{\partial \htH }{\partial \htx_2} \,,
    \label{Eq:shallow3}
\end{eqnarray}
where $g_c = g \,\cos(\alpha)$ and $g_s = g\, \sin(\alpha)$.
We neglect the shear stresses at the bottom and surface, i.e., $\tau_{s1}=\tau_{b1}=\tau_{s2} =\tau_{b2}=0$, and use the Bingham model of Equation \eqref{Eq:bingham}.

Note that the above system is a second-order problem since the terms $\bar{\sigma}_{11}$, $\bar{\sigma}_{12}$, $\bar{\sigma}_{21}$ and $\bar{\sigma}_{22}$ also include spatial derivatives of $\vel_1$ and $\vel_2$.
So, we introduce the auxiliary tensor $\bEE$ \cite{Dawson2005, Hesthaven2017} that represents the gradient of the velocity, i.e.,
\begin{eqnarray}
     \bEE &=& \nabla \bvel
     = \left[ \begin{array}{cc}
         \partial \vel_1 / \partial x_1 & \partial \vel_2 /\partial x_1  \\
         \partial \vel_1 / \partial x_2 & \partial \vel_2 / \partial x_2  \\
    \end{array}\right] 
    \,.
    \label{Eq:gradientVel}
\end{eqnarray}

We rewrite Equations \eqref{Eq:shallow1}, \eqref{Eq:shallow2}, \eqref{Eq:shallow3}, and \eqref{Eq:gradientVel} in a system of first-order equations in vectorial form
\begin{eqnarray}
    \partialt \bfU 
        +  \nabla \cdot  \bfF
        &=& \nabla \cdot  \bfQ +  \bfS  
    \label{Eq:shallowR2}
    \\
    \bfE &=& \nabla \cdot \bfG
    \,,
    \label{Eq:shallowR3}
\end{eqnarray}
where we vectorize the auxiliary tensor $\bEE$ as $\bfE = \vect{ \bEE } = [\EE_{11}\,,  \EE_{21}\,, \EE_{12}\,, \EE_{22}]^\top$, and define the following matrices:
\begin{equation}
    \bfU = \left[
    \begin{array}{c} 
    h \\ h\,\vel_1 \\ 
      h\,\vel_2
     \end{array}\right]\,,
     \label{Eq:Uterm}
\end{equation}
\begin{equation}
    \bfS = \left[
    \begin{array}{c} 
    0 \\ -g_s\, \hth - g_c\,\hth\,  \frac{\partial \htH }{\partial \htx_1} \\ 
     - g_c\,\hth\, \frac{\partial \htH }{\partial \htx_2}
     \end{array}\right]\,,
     \label{Eq:Sterm}
\end{equation}
\begin{eqnarray}
\bfF = \left[ \begin{array}{cc} 
    h\,\vel_1 & h\,\vel_2\\
    h\,\vel_1^2 +1/2 g_c\,\hth^2 & h\,\vel_1 \, \vel_2\\
    h\,\vel_1 \, \vel_2 & h\,\vel_2^2 +1/2 g_c\,\hth^2 
    \end{array}\right] =[\bfF_1\,, \bfF_2]   \,,
    \label{Eq:Fterm}
\end{eqnarray}
\begin{eqnarray}
\bfQ = \left[ \begin{array}{cc} 
    0 & 0\\
    h  \sigma_{11} /\rho & h  \sigma_{12} /\rho\\
    h  \sigma_{21} /\rho & h  \sigma_{22} /\rho 
    \end{array}\right] =[\bfQ_1\,, \bfQ_2]   \,,
    \label{Eq:Qterm}
\end{eqnarray}
and
\begin{eqnarray}
\bfG = \left[ \begin{array}{cc} 
    \vel_1 & 0\\
    \vel_2 & 0\\
    0 & \vel_1 \\
    0 & \vel_2 
    \end{array}\right] =[\bfG_1\,, \bfG_2]   \,.
    \label{Eq:Gterm}
\end{eqnarray}
In this notation. the terms $\nabla \cdot  \bfF$ are $\nabla \cdot  \bfF = \partial\bfF_1/ \partial x_1 + \partial\bfF_2/ \partial x_2$ where $\bfF_1$ and $\bfF_2$ are column vectors.
Note that $\sigma_{ij}$ is a function of $\EE_{ij}$ instead of $\partial \vel_i / \partial x_j$.

\section{Regularized Bingham model}
\label{Sec:Regularization}

In this section, we regularize the Bingham model and express the stress in terms of $\bEE$ instead of $\nabla \bvel$.
First, we rewrite Equation \eqref{Eq:Dtensor} as
\begin{eqnarray}
    \bfD = \frac{1}{2} \left( \bEE + \bEE^{\top} \right)\,,
    \label{Eq:DtensorE}
\end{eqnarray}
With this definition, we split the stress $\bfsigma$ of Equation \eqref{Eq:bingham} in two terms, the Newtonian $\bfsigmaOne$ and Bingham  $\bfsigmaTwo$ contributions
\begin{equation}
    \bfsigma = 
     \bfsigmaOne    +     \bfsigmaTwo\,,
     \label{Eq:binghamR}
\end{equation}
where
\begin{eqnarray}
    \bfsigmaOne=2\eta\, \left( \bfD + \trace (\bfD)\,\bfI \right)\,,
    \label{Eq:sigma1}
\end{eqnarray}
\begin{eqnarray}
    \begin{cases}
        \bfsigmaTwo = \sqrt{2}\sigmao \, \frac{\bfD + \trace (\bfD)\,\bfI}{\normc{\bfD}} & \text{if} \quad \normc{\bfD}\neq 0 \\
       \normc{\bfsigmaTwo} \leq \sqrt{2} \sigmao & \text{if} \quad \normc{\bfD}= 0\,.
    \end{cases}
    \label{Eq:bingham2}
\end{eqnarray}
Note that $\bfsigmaOne= {\bf 0}$ when $\normc{\bfD}= 0$.

\begin{figure}[ht!]
    \centering
    \includegraphics{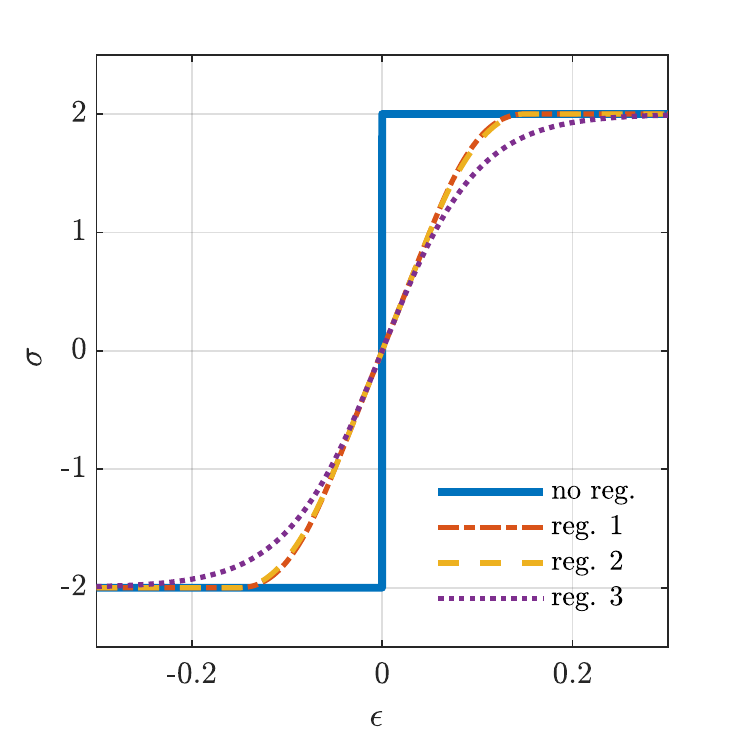}
    \includegraphics{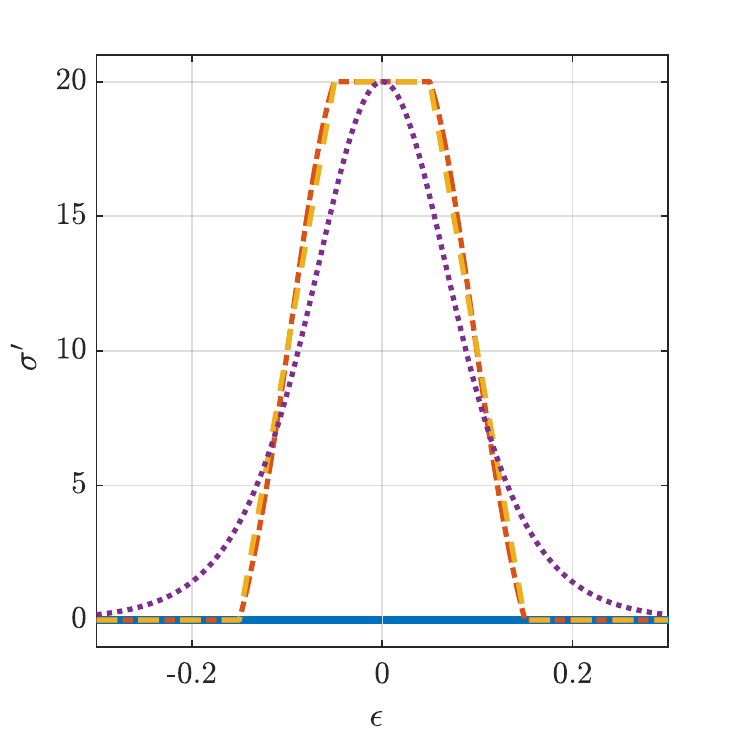}
    \caption{Regularization functions for $\sigmao=1$, $\gamma=10$, and $\beta=1$. }
    \label{Fig:regFunctions}
\end{figure}
\begin{figure}[ht!]
    \centering
    \includegraphics{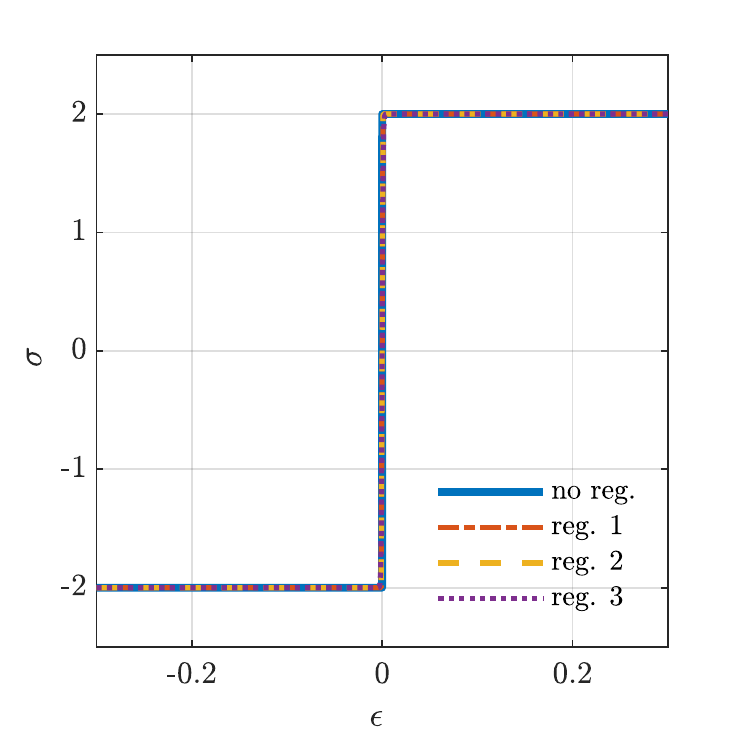}
    \includegraphics{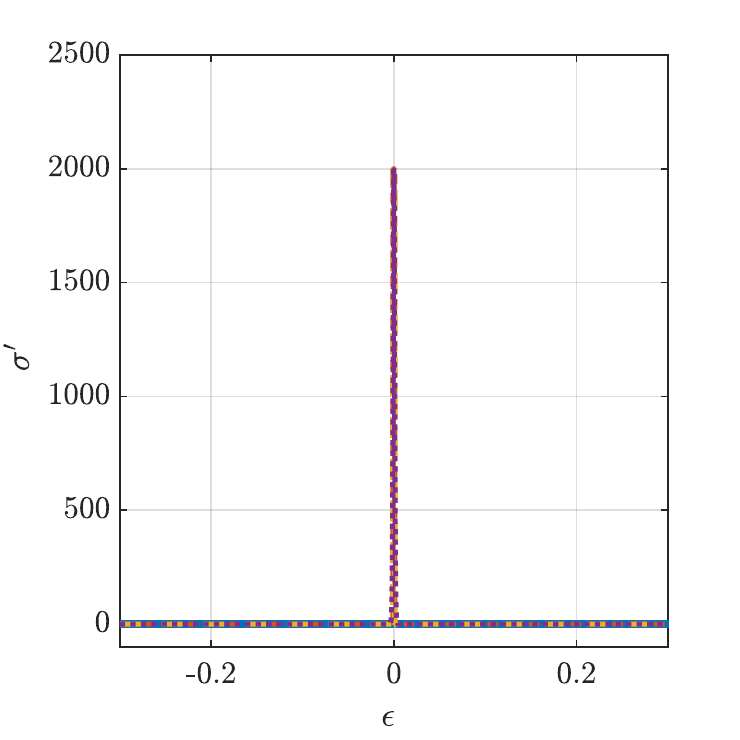}
    \caption{Regularization functions for $\sigmao=1$, $\gamma=10^3$, and $\beta=10^3$. }
    \label{Fig:regFunctionsL}
\end{figure}

For one-dimensional problems, we have $\bfD=\EE=\partial \vel/\partial x$, $\bfD + \trace (\bfD)\,\bfI = 2\,\EE$, $\normc{\bfD}=\sqrt{2}\,|\EE|$, so Equation \eqref{Eq:sigma1} and \eqref{Eq:bingham2} reduce to
\begin{eqnarray}
    \bfsigmaOne=4\eta\, \EE \,,
    \label{Eq:sigma11D}
\end{eqnarray}
\begin{eqnarray}
   \begin{cases}
        \bfsigmaTwoz = 2\sigmao \, \signf{\EE} & \text{if} \quad \left| \EE \right| \ne 0 \\
         \bfsigmaTwoz \le \sigmao  & \text{if} \quad  \EE  = 0\,.
    \end{cases}
    \label{Eq:bingham1D}
\end{eqnarray}
Note that the first derivative of the above Equations \eqref{Eq:bingham1D} is not defined at $\EE  = 0$.
Since we plan to use the Newton-Raphson method to solve our equations, we require to compute the first derivatives.
For this reason, we study three regularizing functions, two of \emph{local} type and one global, with continuous first derivatives.
The first alternative  is to write Equation \eqref{Eq:bingham1D} using a local smooth max function \cite{DeLosReyes2012a} as follows
\begin{eqnarray}
    \bfsigmaTwoa  &=& 2\sigmao \, \frac{ \gamma \, \EE }{\max_{\beta} \left( \gamma\left|\EE \right| \,, \sigmao\right) } \,,
    \label{Eq:sigmaTwoA}
\end{eqnarray}
where $\gamma$ and $\beta$ are regularization parameters and
\begin{eqnarray}
    \max_{\beta} \left( x \,, 0\right) =
    \begin{cases}
        x &\quad \text{if} \quad x\ge \frac{1}{2\beta} \\
        \frac{\beta}{2} \left(x  + \frac{1}{2\beta} \right)^2 &\quad \text{if} \quad |x|\le \frac{1}{2\beta} \\
        0 &\quad \text{if} \quad x \le -\frac{1}{2\beta}.
    \end{cases}
    \label{Eq:sigmaTwoAA}
\end{eqnarray}
The second alternative \cite{DeLosReyes2012a} is also local and is given by the function
\begin{eqnarray}
    \bfsigmaTwob =
    \begin{cases}
        2 \sigmao \frac{\EE}{|\EE|} 
        &\quad \text{if} \quad 
        \gamma|\EE|\ge  \sigmao + \frac{1}{2\beta} 
        \\
        \frac{\EE}{|\EE|} \left(2 \sigmao  -\beta \left(\sigmao-\gamma |\EE|+\frac{1}{2\beta} \right)^2 \right) 
        &\quad \text{if} \quad \sigmao - \frac{1}{2\beta} \le \gamma|\EE|\le  \sigmao + \frac{1}{2\beta} 
        \\
        2\gamma\, \EE 
        &\quad \text{if} \quad \gamma|\EE|\le  \sigmao - \frac{1}{2\beta}.
    \end{cases}
    \label{Eq:sigmaTwoB}
\end{eqnarray}
Finally, as a third alternative, we propose to use a $C^\infty$ global regularization function as follows
\begin{eqnarray}
    \bfsigmaTwoc =
    2 \sigmao \tanh \left( \gamma \,\EE \right) \,.
    \label{Eq:sigmaTwoC}
\end{eqnarray}
For comparison, we plot these regularization functions in Figure \ref{Fig:regFunctions} for $\sigmao=1$, $\gamma=10$, and $\beta=1$. 
Note that all these options have the same slope at $\EE=0$.
The first derivatives are piecewise nonlinear, linear, and smooth for the regularization functions 1, 2, and 3 respectively.
However, in practice we use larger values of the regularization parameters $\gamma$ and $\beta$ to approximate better the reality.
As we show in Figure \ref{Fig:regFunctionsL} for $\sigmao=1$, $\gamma=10^3$, and $\beta=10^3$ the regularized functions resemble our non-regularized function (cf. Equation \eqref{Eq:bingham1D}) but with the main difference that their derivatives are well defined at $\EE=0$.
For completeness, these derivatives are detailed in Appendix \ref{Ape:Derivatives}.

\FloatBarrier
\section{Discontinuous Galerkin}
\label{Sec:DG}

In this section we solve the problem  \eqref{Eq:shallowR2} and \eqref{Eq:shallowR3} given some initial conditions at $t=0$ and prescribed Dirichlet boundary conditions $h^p$, $\bvel^p$, and $\bEE^p$ at the boundaries $\Gamma^{Dh}$, $\Gamma^{D\mu}$, $\Gamma^{D\EE}$ respectively.
We group the trial functions $h$ and $\bvel$ and the test functions as follows:  
\begin{eqnarray}
\begin{array}{ccc}
    \bfV = \left[
    \begin{array}{c} 
    h \\ \vel_1 \\ 
      \vel_2
     \end{array}\right]\,,
&
    \bfw_V = \left[
    \begin{array}{c} 
    w_h \\ w_{\vel 1} \\ 
      w_{\vel 2}
     \end{array}\right]\,,
&
    \bfw_ E = 
    \left[
    \begin{array}{c}
         w_{\EE{11}} \\  w_{\EE{21}}  \\  w_{\EE{12}}  \\  w_{\EE{22}} 
    \end{array}
    \right] 
\,,
\end{array}
     \label{Eq:Vterm}
\end{eqnarray}
where $w_h$, $\bfw_{\vel}$, and $\bfw_E$ are the test functions for the height, velocity, and auxiliary variable, respectively.

In the weak formulation of Equations \eqref{Eq:shallowR2}  and \eqref{Eq:shallowR3}, we define the residual functionals for each equation, i.e., $\bfr_V$ and $\bfr_{E}$ respectively. 
The problem consist in finding $h$, 
$\bvel$
and $\bEE $ such that
\begin{subequations}
    \label{Eq:WeakForm1b}
    \begin{align}
        \bfr_V &= \big\langle \bfwu \,, \partialt \bfU 
        +  \nabla \cdot  \bfF
        -  \nabla \cdot  \bfQ
        - \bfS \big\rangle_{\Omega} = \bfzero \\
        \bfr_E &= \big\langle \bfw_E \,, \bfE - \nabla \cdot \bfG \big\rangle_{\Omega} = \bfzero\,,
    \end{align}
\end{subequations}
for all $w_h $, $\bfw_{\vel}$ and $\bfw_E $.
In this notation, $\langle \bfu, \bfw \rangle_{\Omega} = \int_{\Omega} \bfu \cdot \bfw \, \da$ is the $L^2$ inner product over the set $\Omega$ and $\da$ is the surface differential element.
Note that similarly, $\langle \bfu, \bfw \rangle_{\Gamma} = \int_{\Gamma} \bfu \cdot \bfw \, \dl$ is the $L^2$ inner product over the set $\Gamma$ and $\dl$ is the line differential element.

In the discontinuous Galerkin method, the domain $\Omega$ is discretized in $\nel$ elements such $\Omega^h = \cup_e^{\nel}\Omega_e^h$.
We assume a local solution within each element, and consequently, we enforce no continuity between element boundaries.
We integrate over each element, so System \eqref{Eq:WeakForm1b} is expressed as follows
\begin{subequations}
\label{Eq:WeakForm2}
    \begin{align}
        \bfr_V &= \sum_e \bfr_{Ve} = 
         \sum_e \big\langle  \bfwu \,, \partialt\bfU 
         \big\rangle_{\Omega_e^h}
        + \big\langle \bfwu \,, \nabla \cdot  \bfF \big\rangle_{\Omega_e^h}
        - \big\langle \bfwu \,, \nabla \cdot  \bfQ \big\rangle_{\Omega_e^h}
        -\big\langle \bfwu \,, \bfS \big\rangle_{\Omega_e^h} = \bfzero \\
        \bfr_{E} &= \sum_e \bfr_{E e} = \sum_e \big\langle \bfwe \,, \bfE \big\rangle_{\Omega_e^h} -\big\langle \bfwe \,, \nabla \cdot \bfG \big\rangle_{\Omega_e^h} = \bfzero
        \,.
        \end{align}
\end{subequations}
Since both, the test and trial functions are discontinuous at the element's boundaries, Equations \eqref{Eq:WeakForm2} result in $n_{el}$ local statements
\begin{subequations}
\label{Eq:WeakForm3}
    \begin{align}
        \bfr_{V e} &= \big\langle  \bfwu \,, \partialt \bfU 
         \big\rangle_{\Omega_e^h}
        + \big\langle \bfwu \,, \nabla \cdot  \bfF \big\rangle_{\Omega_e^h}
        - \big\langle \bfwu \,, \nabla \cdot  \bfQ \big\rangle_{\Omega_e^h}
        -\big\langle \bfwu \,, \bfS \big\rangle_{\Omega_e^h} = \bfzero
        \\
        \bfr_{E e} &= \big\langle \bfwe \,, \bfE \big\rangle_{\Omega_e^h} -\big\langle \bfwe \,,\nabla \cdot \bfG \big\rangle_{\Omega_e^h} = \bfzero
        \,,
    \end{align}
\end{subequations}
where $e=1,2,\dots,n_{el}$.

We integrate by parts the flux integrals, for example:
\begin{eqnarray}
    \big\langle \bfwu \,, \bfnabla \cdot  \bfG \big\rangle_{\Omega_e^h} =
    \big\langle \bfwu \,,  \bfG \cdot \bfne \big\rangle_{\Gamma_e^h} 
    -\big\langle  \nabla \bfwu \,, \bfG \big\rangle_{\Omega_e^h} \,,
\end{eqnarray}
where $\bfn_e$ is the unit-normal vector at the boundary of the element $\Gamma_e^h$ pointing outward.
In this notation, $\bfG \cdot \bfne = \bfG_1 \, \bfne\cdot\bfe_1 + \bfG_2 \, \bfne\cdot\bfe_2$ and $\big\langle  \nabla \bfwu \,, \bfG \big\rangle_{\Omega_e^h}= \big\langle  \partial\bfwu/\partial x_1 \,, \bfG_1 \big\rangle_{\Omega_e^h} + \big\langle  \partial\bfwu/\partial x_2 \,, \bfG_2 \big\rangle_{\Omega_e^h}$.
Note that the term $\big\langle \bfwu \,, \bfG \cdot \bfne \big\rangle_{\Gamma_e^h}$ is not well defined  since the trial and test functions are discontinuous at the boundary $\Gamma_e^h$.
To compute this term numerically, first we define the \emph{interior} $\bfVi$ and \emph{exterior} $\bfVe$ values of $\bfV$ at the boundary $\bfx_o \in \Gamma_e^h$ of the element $e$ as follows:
\begin{eqnarray}
    \bfVi(\bfx_o) &=& \lim_{\epsilon\rightarrow 0 \,, \epsilon<0} \bfV(\bfx_o + \epsilon \, \bfne)  \\
    \bfVe(\bfx_o) &=& \lim_{\epsilon\rightarrow 0 \,, \epsilon>0} \bfV(\bfx_o + \epsilon \, \bfne) \,.
\end{eqnarray}
 Similarly, we define the \emph{interior} $\bfwi_V$ 
 \begin{eqnarray}
    \bfwi_V(\bfx_o) &=& \lim_{\epsilon\rightarrow 0 \,, \epsilon<0} \bfw(\bfx_o + \epsilon \, \bfne)  
\end{eqnarray}
and the flux term as follows
\begin{eqnarray}
    \big\langle   \bfwu \,, \bfG(\bfV) \cdot \bfne \big\rangle_{\Gamma_e^h} = \big\langle \bfwi_V \,, \, \nbfG(\bfVi, \bfVe) \, \cdot \bfne\ \big\rangle_{\Gamma_e^h} \,,
\end{eqnarray}
where a \emph{numerical} flux $\nbfG(\bfVi, \bfVe)$ is defined as a function of $\bfUi$ and $\bfUe$.
For example, we can consider the central flux for the term $\nbfG=1/2(\bfGi+ \bfGe)$ where $\bfGi=\bfG(\bfVi)$ and $\bfGe=\bfG(\bfVe)$.
Similarly, we use the central flux for $\nbfQ=1/2(\bfQi+ \bfQe)$ where $\bfQi=\bfQ(\bfVi, \bfEi)$ and $\bfQe=\bfQ(\bfVe, \bfEe)$.
The central flux works well for these terms specifically because there is no preferred direction of propagation.

However, to obtain stable solutions for nonlinear hyperbolic equations, we consider the numerical fluxes as a Riemann problem \cite{Hesthaven2017} because discontinuities are allowed across element boundaries.
The Riemann problem for one-dimensional  conservation laws finds the intercell flux $\bfA(\bfU(x=0),t)$ that solves the following equations.
\begin{subequations}
    \begin{align}
        \partialt \bfU +\nabla\cdot\bfA(\bfU) &= 0 \\
        \bfU(x,0) &= \bfU_L \quad \textrm{if} \quad x\leq0 \\
        \bfU(x,0) &= \bfU_R \quad \textrm{if} \quad x>0\,.
    \end{align}
\end{subequations}
The exact solution to the Riemann problem is relatively straightforward for scalar problems and linear systems. However, the solution for general nonlinear systems is problematic, computationally expensive, and at times impossible \cite{Hesthaven2017}. 
For these reasons, many approximate Riemann solvers have been developed.
Harten, Lax, and van Leer proposed the following HLL Riemann solver 
\begin{eqnarray}
    \nbfA(\bfUi, \bfUe) = \begin{dcases}
        \bfA(\bfUi) \quad \textrm{if} \quad S_L\ge 0
        \\
        \frac{S_R\bfA(\bfUi) -S_L\bfA(\bfUe) + S_LS_R(\bfUe-\bfUi)}{S_R-S_L} \quad \textrm{if} \quad S_L< 0<S_R
        \\
        \bfA(\bfUe) \quad \textrm{if} \quad S_R\le 0
    \end{dcases} \,,
    \label{Eq:HLLSolver}
\end{eqnarray}
where the wave speeds for the one-dimensional shallow-water equations are given by
\begin{subequations}
    \begin{align}
        S_L &= \min\left( \veli -\sqrt{g\hi}, \vele -\sqrt{g\he}\right) \\
        S_R &= \max\left( \veli +\sqrt{g\hi}, \vele +\sqrt{g\he}\right) \,.
    \end{align}
\end{subequations}
Note that this numerical flux considers the direction of the propagation, which is required for the stability of this term.
To compute the numerical flux $\nbfF$, we use the HLL Riemann solver.


From all the above, Equation \eqref{Eq:WeakForm3} is rewritten as
\begin{subequations}
\label{Eq:WeakForm4}
    \begin{align}
       \bfr_{Ve} &=
         \big\langle \bfwu \,,\partialt \bfU 
        \big\rangle_{\Omega_e^h}
        +\big\langle \bfwi_V \,, \nbfF \, \cdot \bfne \big\rangle_{\Gamma_e^h}
        -\big\langle \nabla \bfwu \,, \bfF \big\rangle_{\Omega_e^h}
        \nonumber \\& \quad
        -\big\langle \bfwi_V \,,  \nbfQ\, \cdot \bfne \big\rangle_{\Gamma_e^h }
        +\big\langle \nabla \bfwu \,, \bfQ \big\rangle_{\Omega_e^h}
        -\big\langle \bfwu \,,\bfS \big\rangle_{\Omega_e^h} = \bfzero
        \\
        \bfr_{Ee} &= \big\langle \bfwe \,, \bfE \big\rangle_{\Omega_e^h} -\big\langle \bfwi_E \,,  \nbfG \big\rangle_{ \Gamma_e^h} 
        +\big\langle \nabla \bfwe \,, \bfG \big\rangle_{ \Omega_e^h} 
         = \bfzero
        \,.
    \end{align}
\end{subequations}

We define the local solution for the element $e$, where $\bfx \in \Omega_e^h$, such that
\begin{subequations}
    \begin{align}
        h(\bfx,\htt) &= \sum_{i=1}^{m_1 +1 } \bfh_{e,i}(t) \, \psi_{h,i}(\bfx) =  \hbfpsi{h}(\bfx)^{\top}  \, \bfh_{e}(t)
        \\
        \vel_j(\bfx,\htt) &= \sum_{i=1}^{m_2 +1} \bvel_{j \, e,i}(t) \, \psi_{u,i}(\bfx) =  \hbfpsi{u}(\bfx)^{\top}  \, \bvel_{j\, e}(t)\,, \quad j=\{1,2\}
        \\
        E_j(\bfx,\htt) &= \sum_{i=1}^{m_3 +1} \bfE_{j\,e,i}(t) \, \psi_{E,i}(\bfx) = \hbfpsi{E}(\bfx)^{\top} \, \bfE_{j\, e}(t)\,, \quad j=\{1,2,3,4\}
        \, \,,
    \end{align}
\end{subequations}
where $\psi_{h,i}(\bfx)$, $\psi_{u,i}(\bfx)$ and $\psi_{E,i}(\bfx)$ are the basis functions defined on $\Omega_e^h$.
Note that $h$, $\vel_i$, and $\EE_{ij}$ are defined with polynomials of order $m_1$, $m_2$, and $m_3$ respectively.
For the computations in the examples, we use normalized Legendre polynomial for the basis functions \cite{Hesthaven2017}.
With these definitions, we can easily construct the vectors
\begin{eqnarray}
    \bfV(\bfx, t) = \left[
    \begin{array}{c}
         h \\ \vel_1 \\ \vel_2
    \end{array}
    \right] = 
    \left[
    \begin{array}{ccc}
         \hbfpsi{h}(\bfx)^{\top} & \bfzero_u & \bfzero_u\\ 
         \bfzero_h & \hbfpsi{u}(\bfx)^{\top}  & \bfzero_u \\ 
         \bfzero_h & \bfzero_u & \hbfpsi{u}(\bfx)^{\top}
    \end{array}
    \right] \, 
    \left[
    \begin{array}{c}
         \bfh_e(t) \\ \bvel_{1e}(t) \\ \bvel_{2e}(t)
    \end{array}
    \right]
    = \bfNV (\bfx) \, \bfV_e(t)
    \,,
\end{eqnarray}
\begin{eqnarray}
    \bfE(\bfx, t) = 
    \left[
    \begin{array}{c}
         \EE_{11} \\ \EE_{21} \\ \EE_{12} \\ \EE_{22}
    \end{array}
    \right] = 
    \left[
    \begin{array}{cccc}
         \hbfpsi{E}(\bfx)^{\top} & \bfzero_E & \bfzero_E & \bfzero_E\\ 
         \bfzero_E & \hbfpsi{E}(\bfx)^{\top} &  \bfzero_E & \bfzero_E\\ 
         \bfzero_E & \bfzero_E & \hbfpsi{E}(\bfx)^{\top} & \bfzero_E \\ 
         \bfzero_E & \bfzero_E & \bfzero_E & \hbfpsi{E}(\bfx)^{\top} 
    \end{array}
    \right] \, 
    \left[
    \begin{array}{c}
         \bfE_{1e}(t) \\ \bfE_{2e}(t) \\ \bfE_{3e}(t) \\ \bfE_{4e}(t)
    \end{array}
    \right]
    =\bfNE (\bfx) \,\bfE_e(t)    
    \,,
\end{eqnarray}
where $\bfNV (\bfx)$ and $\bfNE (\bfx)$ are matrices of sizes $3\times(m_1+2m_2+3)$ and $4\times(4m_3+4)$, and $\bfzero_h$, $\bfzero_u$, and $\bfzero_E$ are zero vectors of sizes $1\times(m_1+1)$, $1\times(m_2+1)$, and $1\times(m_3+1)$ respectively.
An important detail is that the bottom height $H$ that appears in the source term for the momentum equation (cf. Equation \eqref{Eq:Sterm}) is given and should be projected to the finite element space for consistency. 
We use polynomials of order $m_0$ to define $H$.
Actually, for conservation and good behaviour of the numerical scheme the height and the bottom height should be discretized in the same way \cite{Bermudez1994}.
Similarly, since we use a Galerkin formulation, the local test functions for the element $e$ are
\begin{subequations}
    \begin{eqnarray}
    \bfw_V(\bfx) = \left[
    \begin{array}{c}
        w_ h \\ w_{\vel 1} \\  w_{\vel 2}
    \end{array}
    \right] = 
     \bfNV (\bfx) \, \bfw_{V\,e}
    \,,
\end{eqnarray}
\begin{eqnarray}
    \bfw_ E(\bfx) = 
    \left[
    \begin{array}{c}
         w_{\EE{11}} \\  w_{\EE{21}}  \\  w_{\EE{12}}  \\  w_{\EE{22}} 
    \end{array}
    \right] 
    =\bfNE (\bfx) \, \bfw_{E\,e}    
    \,,
\end{eqnarray}
\end{subequations}
With these definitions, we express the scheme at each element as follows
\begin{subequations}
    \begin{align}
        \bfr_{Ve} &= \bfwu_{e}^{\top} \left( \, \big\langle \bfNV 
        \, ,  \partialt \bfU \big\rangle_{\Omega_e^h} \,
          + \big\langle \bfNV^- \,, \nbfF \, \cdot \bfne \big\rangle_{\Gamma_e^h} 
         - \big\langle \bfBXV \,, \bfF_1 \big\rangle_{\Omega_e^h}
         - \big\langle \bfBYV \,, \bfF_2 \big\rangle_{\Omega_e^h}
         \right. \nonumber \\& \quad 
         \left.
         -\big\langle \bfNV^-  \,,  \nbfQ \, \cdot \bfne \big\rangle_{\Gamma_e^h}
        +\big\langle \bfBXV \,, \bfQ_1 \big\rangle_{\Omega_e^h}
        +\big\langle \bfBYV \,, \bfQ_2 \big\rangle_{\Omega_e^h}
         -\big\langle \bfNV \,,\bfS \big\rangle_{\Omega_e^h}  \right)
         = \bfwu_{e}^{\top} \bfR_{Ve}= \bfzero 
         \label{Eq:WeakForm51}
         \\
         \bfr_{Ee} &= \bfwe_{e}^{\top} \left( \, \big\langle \bfNE
        \, , \bfE \big\rangle_{\Omega_e^h} \,
         -\big\langle \bfNE^-  \,,  \nbfG \, \cdot \bfne \big\rangle_{\Gamma_e^h}
        +\big\langle \bfBXE \,, \bfG_1 \big\rangle_{\Omega_e^h}
        +\big\langle \bfBYE \,, \bfG_2 \big\rangle_{\Omega_e^h}
        \right) = \bfwe_{e}^{\top} \bfR_{Ee} = \bfzero 
         \,,
         \label{Eq:WeakForm52}
    \end{align}
\end{subequations}
where $\bfBXV = \partial \bfNV / \partial x_1$, $\bfBYV = \partial \bfNV / \partial x_2$, $\bfBXE = \partial \bfNE / \partial x_1$, and $\bfBYE = \partial \bfNE / \partial x_2$.

\subsection{Implicit time integration}

We discretize in time and define $\bfV^{(k)}=\bfV(t^k)$.
The implicit scheme uses $\partialt \bfU(t^{(k)}) =  (\bfU^{(k)} - \bfU^{(k-1)} )/ \Delta t$ and thus the element residuals at time $t^{(k)}$ are given by
\begin{subequations}
    \label{Eq:WeakForm7}
    \begin{align}
        \bfR_{Ve}^{(k)} &=  \, \frac{1}{\Delta t}\big\langle \bfNV 
        \, ,  \bfU^{(k)} - \bfU^{(k-1)} \big\rangle_{\Omega_e^h} \,
          + \big\langle \bfNV ^- \,, \nbfF^{(k)} \, \cdot \bfne \big\rangle_{\Gamma_e^h} 
         - \big\langle \bfBXV \,, \bfF_1^{(k)} \big\rangle_{\Omega_e^h}
         - \big\langle \bfBYV \,, \bfF_2^{(k)} \big\rangle_{\Omega_e^h}
          \nonumber \\& \quad 
         -\big\langle \bfNV ^-  \,,  \nbfQ^{(k)} \, \cdot \bfne \big\rangle_{\Gamma_e^h}
        +\big\langle \bfBXV \,, \bfQ_1^{(k)} \big\rangle_{\Omega_e^h}
        +\big\langle \bfBYV \,, \bfQ_2^{(k)} \big\rangle_{\Omega_e^h}
         -\big\langle \bfNV  \,,\bfS^{(k)} \big\rangle_{\Omega_e^h}  
         = \bfzero 
         \\
         \bfR_{Ee}^{(k)} &=  \, \big\langle \bfNE 
        \, , \bfE^{(k)} \big\rangle_{\Omega_e^h} \,
         -\big\langle \bfNE^-  \,,  \nbfG^{(k)} \, \cdot \bfne \big\rangle_{\Gamma_e^h}
        +\big\langle \bfBXE \,, \bfG_1^{(k)} \big\rangle_{\Omega_e^h}
        +\big\langle \bfBYE \,, \bfG_2^{(k)} \big\rangle_{\Omega_e^h}
          = \bfzero
         \,,
    \end{align}
\end{subequations}
To use Newton-Raphson's method to solve \eqref{Eq:WeakForm7} for $\bfV^{(k)}$ and $\bfE^{(k)}$, we compute the tangent terms ${\partial \bfR_{Ve}^{(k)}}/{\partial \bfV_e}$, ${\partial \bfR_{Ve}^{(k)}}/{\partial \bfE_e}$, ${\partial \bfR_{Ee}^{(k)}}/{\partial \bfV_e}$, and ${\partial \bfR_{Ee}^{(k)}}/{\partial \bfE_e}$.
For completeness, these derivatives are detailed in Appendix \ref{Ape:Derivatives}.
We solve for $\bfV^{(k)}$, and $\bfE^{(k)}$ via Newton-Raphson's method, whereupon we update the current solution guess $\bfV^I(\bfx,\,t^{(k)})$ and $\bfE^I(\bfx,\,t^{(k)})$ to $\bfV^{I+1}(\bfx,\,t^{(k)})=\bfV^I(\bfx,\,t^{(k)})+\Delta \bfV$
and $\bfE^{I+1}(\bfx,\,t^{(k)})=\bfE^I(\bfx,\,t^{(k)})+\Delta \bfE$ 
where $\Delta \bfV$ and $\Delta \bfE$ are the solution to the linear problem: 
\begin{eqnarray}
    \left[ \begin{array}{cc} 
    \partial \bfR_{Ve}^{(k)}/\partial \bfV & 
    \partial \bfR_{Ve}^{(k)}/\partial \bfE \\
    \partial \bfR_{Ee}^{(k)}/\partial \bfV & 
    \partial \bfR_{Ee}^{(k)}/\partial \bfE \\
    \end{array}\right] 
    \left[ \begin{array}{c} 
    \Delta\bfV\\
    \Delta \bfE 
    \end{array}\right]
    = -
    \left[ \begin{array}{c} 
    \bfR_{Ve}^{(k)}\\
    \bfR_{Ee}^{(k)}
    \end{array}\right].
    \label{Eq:Newton}
\end{eqnarray}

\subsection{Boundary conditions}
We impose strongly the boundary conditions where we define the external values of our variables at the boundary accordingly. 
Thus, for the Dirichlet boundary condition, we simply define
\begin{equation}
    \bfVe(\bfx_o) = \bfV^p \,,
\end{equation}
where $\bfV^p$ is the prescribed value at the boundary and $\bfx_o \in \Gamma$. 
Similarly, for the homogeneous Neumann boundary condition, we have
\begin{equation}
    \bfVe(\bfx_o) = \bfVi(\bfx_o) \,,
\end{equation}
where $\nabla \bfV \cdot \bfn(\bfx_o)  = 0$ is imposed at the boundary point $\bfx_o \in \Gamma$.
We impose the boundary conditions similarly for the auxiliary variable vector $\bfE$.


\FloatBarrier
\section{Examples}
In this section, we solve three example problems to show the effectiveness and potential of our method.
The first two examples test the well-balanced properties of our approach.
In all these examples, we also validate the accuracy and computational cost of the different regularization functions proposed in Section \ref{Sec:Regularization}.
We set the maximum number of Newton-Raphson iterations per time step as 10 for all the problems.
As mentioned before, we use normalized Legendre polynomials for our basis functions, and we use the HLL Riemann solver to compute the numerical flux $\nbfF$.
For all the regularization approaches, if the material strain rate magnitude is $|\EE|>\sigmao/\gamma$ in a portion of the domain, we consider it an active region.
Thus, the percentage of the active region in the domain is computed as
\begin{eqnarray}
    \textrm{active\%} = \frac{100 \%}{|\Omega|}\int_{\Omega} \Heav(|\EE|-\sigmao/\gamma) \, \da \,,
\end{eqnarray}
where $|\Omega|=\int_{\Omega} \da$, and $\Heav$ is the Heaviside function, i.e., $\Heav(x)=1$ if $x\ge0$ and $\Heav(x)=0$ if $x<0$.

\label{Sec:Examples}
\subsection{Constant free surface}
\begin{figure}
    \centering
    \includegraphics{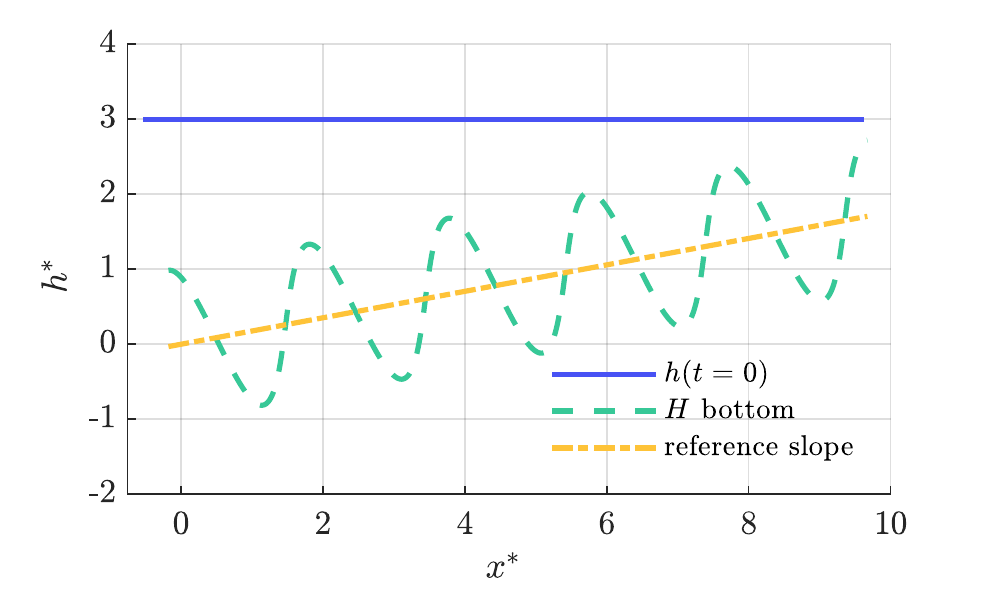}
    \caption{Constant free surface test problem.}
    \label{Fig:test1}
\end{figure}

In this first example, we test if our discontinuous Galerkin approach preserves the initial equilibrium condition, i.e., method is well-balanced.
We consider a domain of length $\Omega=[0,L]$ where $L=10$ with the following initial condition
\begin{subequations}
    \begin{align}
        \vel(x,0) &= 0 \\
        h(x,0) &= \frac{h_0}{\cos(\alpha)} - x \tan(\alpha) - H(x) \\
        H(x) &= \cos(\pi x)\,,
    \end{align}
\end{subequations}
where the slope is $\alpha=\pi/18$ and $h_0=3$, 
cf. Figure \ref{Fig:test1}.
For this problem, we apply only Dirichlet boundary conditions on both sides, i.e., $h(0,t)=h_0/\cos(\alpha)-1$, $h(L,t)=h_0/\cos(\alpha)-L\tan(\alpha)-\cos(\pi L)$,  $\vel(0,t)=\EE(0,t)=\vel(L,t)=\EE(L,t)=0$.
We do not expect any movement in this problem because the fluid is at rest and the free surface is constant and horizontal (perpendicular to the gravity vector).

The solution is computed at time $t=1$ using different polynomial orders $m_0, m_1, m_2, m_3 \in \{1,2,3,4\}$ for the projected bottom height $H$ and our variables $h$, $\vel$, and $\EE$ respectively.
We discretize the domain with $n_{el}\in \{100 ,\, 1000\}$ elements and use the time step $\Delta t=1e-2$.
The viscosity and Bingham yield stress is set for $\eta\in\{0, \, 1\}$ and $\sigmao\in\{0, \, 1\}$ respectively.
We use the regularization approach 1 (cf. Equation \eqref{Eq:sigmaTwoA}) with the regularization parameter $\gamma = \beta= 10^3$ to solve this problem.
The density and gravity are $\rho=1$ and $g=9.81$.
As mentioned, since the initial condition is in equilibrium (satisfies Equations \eqref{Eq:shallowR2} and \eqref{Eq:shallowR3}), it must be preserved if our approach is well-balanced.
For this reason, we compute the error norms L-2 and L-infinity for the height and velocity denoted as \Lth, \Ltvel, \Lh, and \Lvel \,  by comparing  the quadrature points to the initial condition.
Results are summarized in Tables \ref{tab:errorsEx1a} and \ref{tab:errorsEx1b}.

As expected, our approach preserves the initial condition as long as the polynomial order of the basis that describes the bottom is smaller or equal to the polynomial order used for the height, i.e., $m_0 \le m_1$.
However, the error decreases for a finer mesh or a higher order height polynomial $m_1$ even if $m_0>m_1$ (shaded in gray on both tables Tables \ref{tab:errorsEx1a} and \ref{tab:errorsEx1b}).
This behavior makes sense because the height needs to be described at least with the same order polynomial as the bottom to cancel the effect of the bottom source term (cf. Equation \eqref{Eq:Sterm}) in the equations.
Otherwise, for $m_0>m_1$ the equations might not be numerically satisfied, so the heights are slightly different from the initial condition, and the velocity is not exactly zero.
Note that for $m_0 \le m_1$, the height and zero velocity are preserved regardless of any parameter or number of elements.
This preservation of the initial condition is also independent of the viscosity and Bingham yield stress.

\begin{figure}
    \centering
    \includegraphics{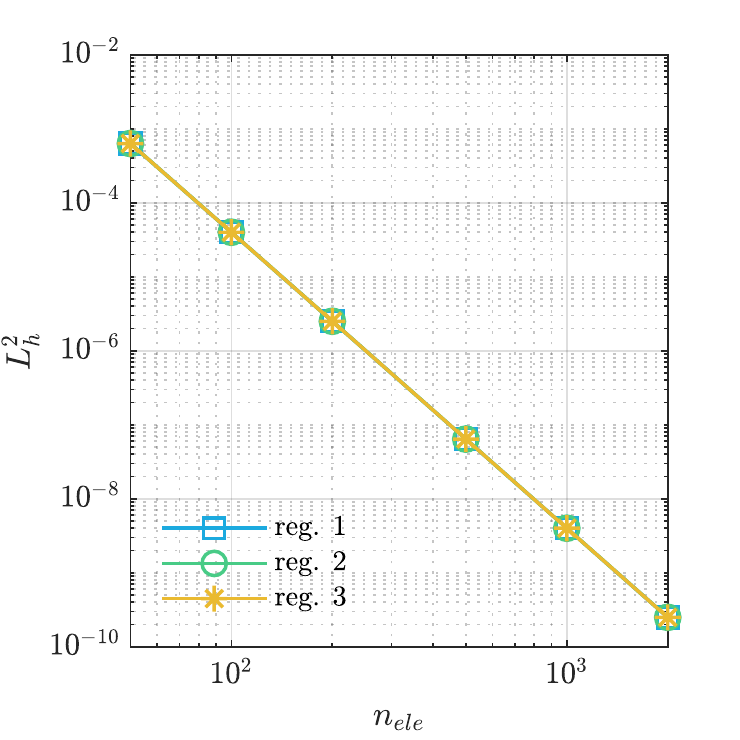}
    \includegraphics{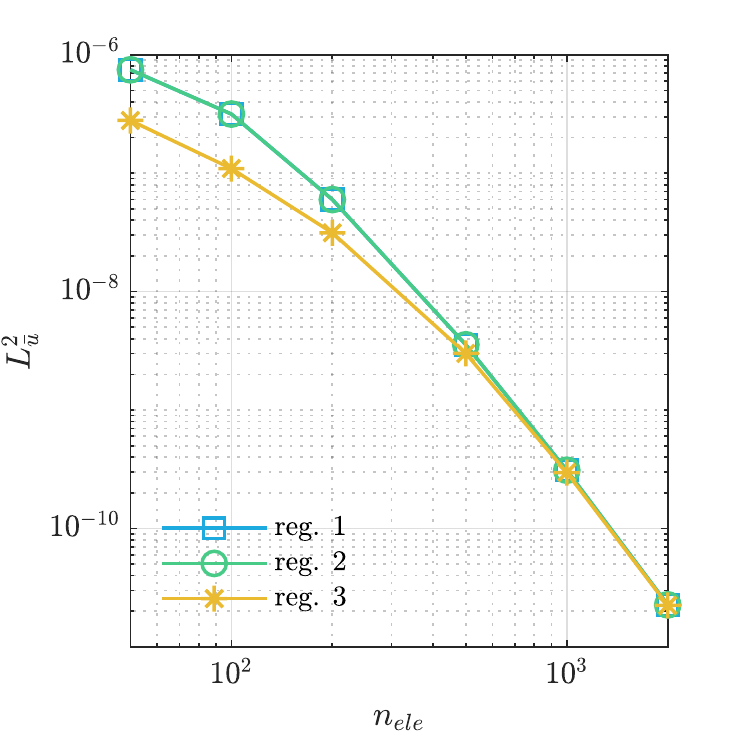}
    \includegraphics{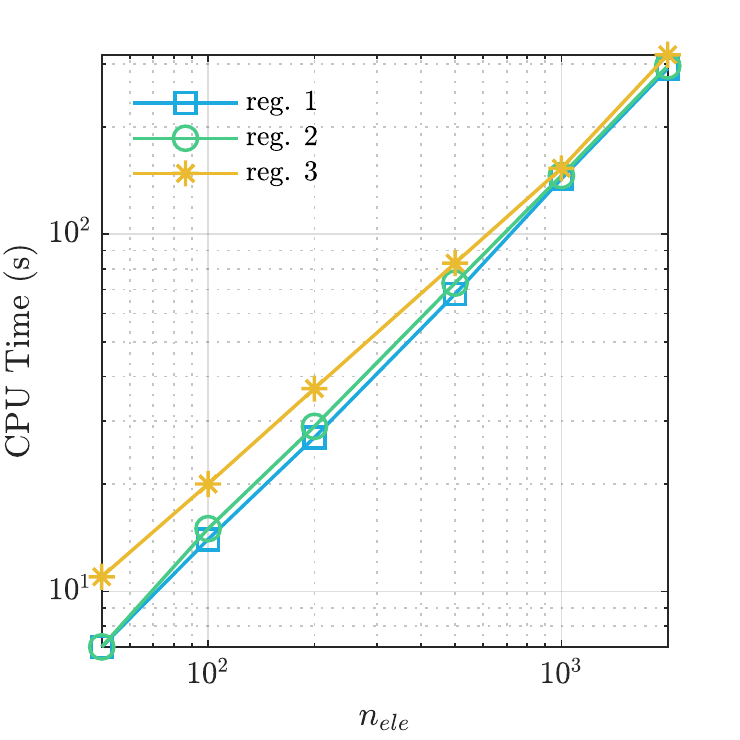}
    \caption{Preservation error  and CPU time for the test problem with constant free surface using $m_0=2$ and $m_1=m_2=m_3=1$ and different regularization approaches.}
    \label{Fig:errorTest1}
\end{figure}

We repeated the same experiments with the other regularizations and obtained the same behavior in the results.
For $m_0 \le m_1$, the height and zero velocity are preserved regardless of any parameter, number of elements, 
and regularization approach.
Interestingly, for $m_0>m_1$, the errors have a subtle difference depending on the regularization approach.
We present these results for the different regularization approaches in Figure \ref{Fig:errorTest1} and Table \ref{tab:errorsEx1c} using $m_0=2$ and $m_1=m_2=m_3=1$.
Note that for this case, the height error is the same for all the approaches.
As expected, for regularization approaches 1, 2, and 3 present quadratic convergence, taking 2-3 Newton-Raphson iterations per time iteration.
Another interesting point is that for a coarse mesh, the third regularization approach has a lower velocity error but requires more computational time than the other methods.
For a finer mesh, all the regularization approaches have the same result accuracy and computational cost for this problem.

\begin{table}[!ht]
    \centering
    \begin{tabular}{ccccccccccc}
         $n_{el}$ &  $m_0$ &  $m_1$ &  $m_2$ &  $m_3$ &  $\eta$ & $\sigmao$ & \Lth & \Ltvel & \Lh & \Lvel \\
\hline
 100 &  1 &  1 &  1 &  1 &     0 &     0 & 9.8E-32 & 0.0E+00 & 8.8E-17 & 0.0E+00 \\ 
      100 &  2 &  2 &  2 &  2 &     0 &     0 & 7.5E-32 & 0.0E+00 & 8.8E-17 & 0.0E+00 \\ 
      100 &  3 &  3 &  3 &  3 &     0 &     0 & 1.8E-17 & 1.9E-16 & 2.7E-09 & 7.8E-09 \\ 
      100 &  4 &  4 &  4 &  4 &     0 &     0 & 1.4E-21 & 1.2E-19 & 2.8E-11 & 2.8E-10 \\ 
     \hline 
      100 &  1 &  1 &  1 &  1 &     1 &     0 & 9.8E-32 & 0.0E+00 & 8.8E-17 & 0.0E+00 \\ 
      100 &  2 &  2 &  2 &  2 &     1 &     0 & 7.5E-32 & 0.0E+00 & 8.8E-17 & 0.0E+00 \\ 
      100 &  3 &  3 &  3 &  3 &     1 &     0 & 8.6E-19 & 1.3E-16 & 2.9E-10 & 5.2E-09 \\ 
      100 &  4 &  4 &  4 &  4 &     1 &     0 & 7.1E-26 & 1.1E-23 & 9.9E-14 & 1.7E-12 \\ 
     \hline 
      100 &  1 &  1 &  1 &  1 &     1 &     1 & 9.8E-32 & 0.0E+00 & 8.8E-17 & 0.0E+00 \\ 
      100 &  2 &  2 &  2 &  2 &     1 &     1 & 7.5E-32 & 0.0E+00 & 8.8E-17 & 0.0E+00 \\ 
      100 &  3 &  3 &  3 &  3 &     1 &     1 & 4.9E-20 & 2.6E-17 & 6.5E-11 & 1.5E-09 \\ 
      100 &  4 &  4 &  4 &  4 &     1 &     1 & 4.2E-31 & 1.4E-28 & 3.1E-16 & 7.8E-15 \\ 
     \hline 
     \rowcolor{Gray}
      100 &  2 &  1 &  1 &  1 &     1 &     1 & 4.0E-05 & 3.2E-07 & 1.2E-03 & 1.4E-04 \\ 
      100 &  2 &  2 &  1 &  1 &     1 &     1 & 1.9E-11 & 1.2E-12 & 1.8E-06 & 2.2E-07 \\ 
      \rowcolor{Gray}
      100 &  2 &  1 &  2 &  1 &     1 &     1 & 4.0E-05 & 1.9E-04 & 1.2E-03 & 4.4E-03 \\ 
      \rowcolor{Gray}
      100 &  2 &  1 &  1 &  2 &     1 &     1 & 4.0E-05 & 8.5E-11 & 1.2E-03 & 1.4E-06 \\ 
     \hline 
     \rowcolor{Gray}
      100 &  3 &  2 &  2 &  2 &     1 &     1 & 1.6E-08 & 8.6E-16 & 2.4E-05 & 7.5E-09 \\ 
      100 &  3 &  3 &  2 &  2 &     1 &     1 & 9.9E-23 & 3.9E-21 & 5.0E-12 & 4.8E-11 \\ 
      \rowcolor{Gray}
      100 &  3 &  2 &  3 &  2 &     1 &     1 & 1.6E-08 & 2.4E-09 & 2.4E-05 & 2.3E-05 \\ 
      \rowcolor{Gray}
      100 &  3 &  2 &  2 &  3 &     1 &     1 & 1.6E-08 & 2.7E-16 & 2.4E-05 & 3.3E-09 \\ 
     \hline 
     \rowcolor{Gray}
      100 &  4 &  3 &  3 &  3 &     1 &     1 & 7.5E-12 & 2.2E-13 & 5.0E-07 & 1.4E-07 \\ 
      100 &  4 &  4 &  3 &  3 &     1 &     1 & 4.8E-18 & 9.6E-20 & 6.5E-10 & 8.8E-11 \\ 
      \rowcolor{Gray}
      100 &  4 &  3 &  4 &  3 &     1 &     1 & 8.5E-12 & 4.6E-11 & 5.5E-07 & 3.3E-06 \\ 
      \rowcolor{Gray}
      100 &  4 &  3 &  3 &  4 &     1 &     1 & 7.5E-12 & 3.0E-20 & 5.0E-07 & 4.4E-11 \\ 
     \hline 
      100 &  1 &  2 &  2 &  2 &     1 &     1 & 6.4E-32 & 0.0E+00 & 8.8E-17 & 0.0E+00 \\ 
      100 &  2 &  3 &  3 &  3 &     1 &     1 & 1.7E-31 & 0.0E+00 & 1.3E-16 & 0.0E+00 \\ 
      100 &  3 &  4 &  4 &  4 &     1 &     1 & 1.1E-31 & 0.0E+00 & 1.8E-16 & 0.0E+00 \\ 
      100 &  1 &  2 &  1 &  2 &     1 &     1 & 6.4E-32 & 0.0E+00 & 8.8E-17 & 0.0E+00 \\ 
      100 &  2 &  3 &  2 &  3 &     1 &     1 & 1.7E-31 & 0.0E+00 & 1.3E-16 & 0.0E+00 \\ 
      100 &  3 &  4 &  3 &  4 &     1 &     1 & 2.8E-27 & 1.4E-25 & 1.1E-14 & 1.3E-13 \\ 
     \hline 
     \hline 
    \end{tabular}
    \caption{Preservation of the initial equilibrium for the test problem with constant free surface using $n_{el}=100$ elements.}
    \label{tab:errorsEx1a}
\end{table}

\begin{table}[!ht]
    \centering
    \begin{tabular}{ccccccccccc}
         $n_{el}$ &  $m_0$ &  $m_1$ &  $m_2$ &  $m_3$ &  $\eta$ & $\sigmao$ & \Lth & \Ltvel & \Lh & \Lvel \\
         \hline
1000 &  1 &  1 &  1 &  1 &     0 &     0 & 9.2E-32 & 0.0E+00 & 8.8E-17 & 0.0E+00 \\ 
     1000 &  2 &  2 &  2 &  2 &     0 &     0 & 9.7E-32 & 0.0E+00 & 1.8E-16 & 0.0E+00 \\ 
     1000 &  3 &  3 &  3 &  3 &     0 &     0 & 2.0E-31 & 0.0E+00 & 1.8E-16 & 0.0E+00 \\ 
     1000 &  4 &  4 &  4 &  4 &     0 &     0 & 1.1E-31 & 0.0E+00 & 1.8E-16 & 0.0E+00 \\ 
     \hline 
     1000 &  1 &  1 &  1 &  1 &     1 &     0 & 9.2E-32 & 0.0E+00 & 8.8E-17 & 0.0E+00 \\ 
     1000 &  2 &  2 &  2 &  2 &     1 &     0 & 9.7E-32 & 0.0E+00 & 1.8E-16 & 0.0E+00 \\ 
     1000 &  3 &  3 &  3 &  3 &     1 &     0 & 2.0E-31 & 0.0E+00 & 1.8E-16 & 0.0E+00 \\ 
     1000 &  4 &  4 &  4 &  4 &     1 &     0 & 1.1E-31 & 0.0E+00 & 1.8E-16 & 0.0E+00 \\ 
     \hline 
     1000 &  1 &  1 &  1 &  1 &     1 &     1 & 9.2E-32 & 0.0E+00 & 8.8E-17 & 0.0E+00 \\ 
     1000 &  2 &  2 &  2 &  2 &     1 &     1 & 9.7E-32 & 0.0E+00 & 1.8E-16 & 0.0E+00 \\ 
     1000 &  3 &  3 &  3 &  3 &     1 &     1 & 2.0E-31 & 0.0E+00 & 1.8E-16 & 0.0E+00 \\ 
     1000 &  4 &  4 &  4 &  4 &     1 &     1 & 1.1E-31 & 0.0E+00 & 1.8E-16 & 0.0E+00 \\ 
     \hline 
     \rowcolor{Gray}
     1000 &  2 &  1 &  1 &  1 &     1 &     1 & 4.0E-09 & 3.1E-10 & 1.2E-05 & 4.2E-06 \\ 
     1000 &  2 &  2 &  1 &  1 &     1 &     1 & 2.0E-19 & 2.3E-20 & 1.9E-10 & 8.3E-11 \\ 
     \rowcolor{Gray}
     1000 &  2 &  1 &  2 &  1 &     1 &     1 & 4.0E-09 & 2.3E-08 & 1.2E-05 & 4.9E-05 \\ 
     \rowcolor{Gray}
     1000 &  2 &  1 &  1 &  2 &     1 &     1 & 4.0E-09 & 8.6E-15 & 1.2E-05 & 1.4E-08 \\ 
     \hline 
     \rowcolor{Gray}
     1000 &  3 &  2 &  2 &  2 &     1 &     1 & 1.6E-14 & 2.7E-24 & 2.4E-08 & 4.0E-13 \\ 
     1000 &  3 &  3 &  2 &  2 &     1 &     1 & 2.0E-31 & 0.0E+00 & 1.8E-16 & 0.0E+00 \\ 
     \rowcolor{Gray}
     1000 &  3 &  2 &  3 &  2 &     1 &     1 & 1.6E-14 & 2.5E-17 & 2.4E-08 & 2.4E-09 \\ 
     \rowcolor{Gray}
     1000 &  3 &  2 &  2 &  3 &     1 &     1 & 1.6E-14 & 2.5E-24 & 2.4E-08 & 3.1E-13 \\ 
     \hline 
     \rowcolor{Gray}
     1000 &  4 &  3 &  3 &  3 &     1 &     1 & 7.5E-20 & 3.8E-21 & 5.0E-11 & 2.2E-11 \\ 
     1000 &  4 &  4 &  3 &  3 &     1 &     1 & 1.1E-31 & 0.0E+00 & 1.8E-16 & 0.0E+00 \\ 
     \rowcolor{Gray}
     1000 &  4 &  3 &  4 &  3 &     1 &     1 & 7.5E-20 & 4.2E-19 & 5.0E-11 & 3.1E-10 \\ 
     \rowcolor{Gray}
     1000 &  4 &  3 &  3 &  4 &     1 &     1 & 7.5E-20 & 2.9E-30 & 5.0E-11 & 4.5E-16 \\ 
     \hline 
     1000 &  1 &  2 &  2 &  2 &     1 &     1 & 7.9E-32 & 0.0E+00 & 1.3E-16 & 0.0E+00 \\ 
     1000 &  2 &  3 &  3 &  3 &     1 &     1 & 2.0E-31 & 0.0E+00 & 2.2E-16 & 0.0E+00 \\ 
     1000 &  3 &  4 &  4 &  4 &     1 &     1 & 1.0E-31 & 0.0E+00 & 1.8E-16 & 0.0E+00 \\ 
     1000 &  1 &  2 &  1 &  2 &     1 &     1 & 7.9E-32 & 0.0E+00 & 1.3E-16 & 0.0E+00 \\ 
     1000 &  2 &  3 &  2 &  3 &     1 &     1 & 2.0E-31 & 0.0E+00 & 2.2E-16 & 0.0E+00 \\ 
     1000 &  3 &  4 &  3 &  4 &     1 &     1 & 1.0E-31 & 0.0E+00 & 1.8E-16 & 0.0E+00 \\ 
     \hline  
     \hline 
    \end{tabular}
    \caption{Preservation of the initial equilibrium for the test problem with constant free surface using $n_{el}=1000$ elements.}
    \label{tab:errorsEx1b}
\end{table}

\begin{table}[!ht]
    \centering
    \begin{tabular}{cccccccc}
         Reg. &  $n_{ele}$ & \Lth & \Ltvel & \Lh & \Lvel & NRIt/Tit & CPUTime\\
         approach &   &   &   &   &   &  & (s) \\
     \hline 
        1 &   50 & 6.3E-04 & 7.5E-07 & 4.6E-03 & 2.2E-04 &  2.170 & 8.0E+00\\ 
        1 &  100 & 4.0E-05 & 3.2E-07 & 1.2E-03 & 1.4E-04 &  2.090 & 1.6E+01\\ 
        1 &  200 & 2.5E-06 & 6.0E-08 & 3.1E-04 & 6.0E-05 &  2.000 & 3.0E+01\\ 
        1 &  500 & 6.4E-08 & 3.6E-09 & 4.9E-05 & 1.4E-05 &  1.990 & 7.5E+01\\ 
        1 & 1000 & 4.0E-09 & 3.1E-10 & 1.2E-05 & 4.2E-06 &  1.980 & 1.6E+02\\ 
        1 & 2000 & 2.5E-10 & 2.3E-11 & 3.1E-06 & 1.1E-06 &  1.980 & 3.2E+02\\ 
     \hline 
        2 &   50 & 6.3E-04 & 7.5E-07 & 4.6E-03 & 2.2E-04 &  2.170 & 8.0E+00\\ 
        2 &  100 & 4.0E-05 & 3.2E-07 & 1.2E-03 & 1.4E-04 &  2.090 & 1.8E+01\\ 
        2 &  200 & 2.5E-06 & 6.0E-08 & 3.1E-04 & 6.0E-05 &  2.000 & 3.1E+01\\ 
        2 &  500 & 6.4E-08 & 3.6E-09 & 4.9E-05 & 1.4E-05 &  1.990 & 7.9E+01\\ 
        2 & 1000 & 4.0E-09 & 3.1E-10 & 1.2E-05 & 4.2E-06 &  1.980 & 1.6E+02\\ 
        2 & 2000 & 2.5E-10 & 2.3E-11 & 3.1E-06 & 1.1E-06 &  1.980 & 3.3E+02\\ 
     \hline 
        3 &   50 & 6.3E-04 & 2.8E-07 & 4.6E-03 & 1.3E-04 &  3.050 & 1.1E+01\\ 
        3 &  100 & 4.0E-05 & 1.1E-07 & 1.2E-03 & 8.3E-05 &  2.770 & 2.1E+01\\ 
        3 &  200 & 2.5E-06 & 3.1E-08 & 3.0E-04 & 4.4E-05 &  2.540 & 4.2E+01\\ 
        3 &  500 & 6.4E-08 & 3.0E-09 & 4.9E-05 & 1.3E-05 &  2.230 & 9.0E+01\\ 
        3 & 1000 & 4.0E-09 & 3.0E-10 & 1.2E-05 & 4.1E-06 &  2.070 & 1.6E+02\\ 
        3 & 2000 & 2.5E-10 & 2.2E-11 & 3.1E-06 & 1.1E-06 &  2.000 & 3.3E+02\\ 
     \hline 
    \end{tabular}
    \caption{Preservation error  and CPU time for the test problem with constant free surface using $m_0=2$ and $m_1=m_2=m_3=1$ and different regularization approaches.}
    \label{tab:errorsEx1c}
\end{table}

\FloatBarrier

\subsection{Parallel free surface to the reference slope}
\begin{figure}
    \centering
    \includegraphics{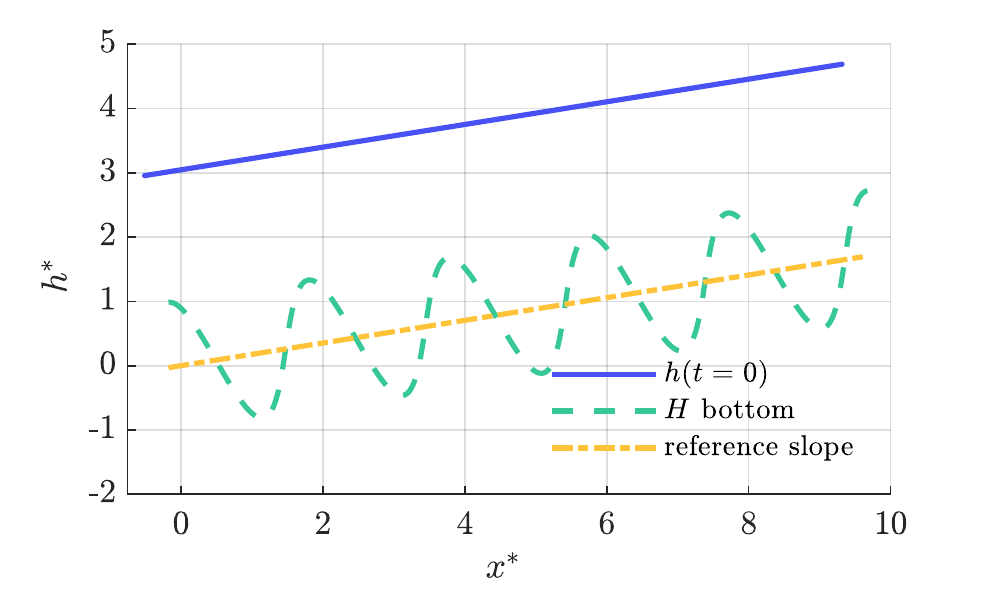}
    \caption{Parallel free surface to the reference slope problem.}
    \label{Fig:test2}
\end{figure}

In the initial condition for this example, the free surface is parallel to the reference slope, and the bottom is variable. 
We test if our discontinuous Galerkin approach preserves the initial equilibrium condition if the material is rigid enough.
To do this, we consider a domain of length $\Omega=[0,L]$ where $L=10$ with the following initial condition
\begin{subequations}
    \begin{align}
        \vel(x,0) &= 0 \\
        h(x,0) &= h_0 - H(x) \\
        H(x) &= \cos(\pi x)\,,
    \end{align}
\end{subequations}
where the slope is $\alpha=\pi/18$ and $h_0=3$, 
cf. Figure \ref{Fig:test2}.
For this problem, we apply the following boundary conditions
\begin{subequations}
    \begin{align}
        \frac{\partial h}{\partial x}(0,t) &= 0, \quad & \frac{\partial \vel}{\partial x}(0,t) &=0, \quad &\frac{\partial \EE}{\partial x}(0,t) &= 0, \\
        \frac{\partial h}{\partial x}(L,t) &= 0, \quad &
        \vel(L,t) &=0, \quad &
        \frac{\partial \EE}{\partial x}(L,t) &= 0 \,.
    \end{align}
\end{subequations}
The density, slope angle, and gravity are $\rho=1$, $\alpha=\pi/18$, and $g=9.81$, respectively.

Since the fluid is at rest, no movement is expected if the material's Bingham yield stress is above certain threshold, i.e., satisfies the following general condition
\begin{eqnarray}
    \left|-\htrho \, g \sin{\alpha} \int_{L/2}^x H(s) \, \textrm{d}s\right| \le \sigmao \sqrt{2} H(x) \quad \forall x \in [0,L]\,.
\end{eqnarray}
For details about this condition, the reader is referred to \cite{Fernandez2014}.
In this specific example, the condition is given by
\begin{eqnarray}
    \sigmao &\ge& \max_{x\in[0,L]} \frac{\htrho \, g\sin(\alpha)}{\sqrt{2}}\frac{\left|  h_0 \left( x-L/2 \right) - \sin(\pi x)/ \pi + \sin(\pi L)/\pi \right|}{h_0 - \cos(\pi x)}\,.
\\
    \sigmao &\ge& \frac{\htrho \, g \,\sin(\alpha) \, h_0  \, L}{2\sqrt{2}(h_0 - 1)} \approx 9.0341\,.
    \label{Eq:conditonTest2}
\end{eqnarray}

The solution is computed at time $t=1$ using polynomial orders $m_0= m_1= m_2= m_3=1$ otherwise stated for the projected bottom height $H$ and our variables $h$, $\vel$, $\EE$ respectively.
We discretize the domain with $n_{el}=100$ and use the time step $\Delta t=1e-6$.
We set the viscosity to $\eta=1$,  and the Bingham yield stress is $\sigmao=9.035$, which is slightly larger than the condition given in Equation \eqref{Eq:conditonTest2}.
Since the initial condition is in equilibrium and the Bingham yield stress is above the threshold, the fluid should preserve the initial condition.
We compute the error norms L-2 and L-infinity for the height and velocity denoted as \Lth, \Ltvel, \Lh, and \Lvel at the quadrature points comparing them to the initial condition.
We solve this problem using the regularization functions described in Equations  \eqref{Eq:sigmaTwoA}, \eqref{Eq:sigmaTwoB}, and \eqref{Eq:sigmaTwoC}.
Results are summarized in Table \ref{tab:errorsEx2}.

\begin{table}[!ht]
    \centering
    \begin{tabular}{cccccccccc}
         Reg. &  $\gamma$ &  $\beta$ &  active &  \Lth & \Ltvel & \Lh & \Lvel & NRIt/TIt & CPUTime \\
          approach &   &   & $\%$ &   &  &  &  &  & (s) \\
     \hline 
     1 &   1E+02 &  1E+02 &     0 & 1.7E-01 & 1.0E-01 & 2.0E-01 & 1.2E-01 &      2 &1.3E+05\\ 
     1 &   1E+02 &  1E+03 &     0 & 1.7E-01 & 1.0E-01 & 2.0E-01 & 1.2E-01 &      2 &1.3E+05\\ 
     1 &   1E+03 &  1E+02 &     0 & 1.8E-03 & 9.7E-04 & 2.1E-02 & 1.3E-02 &      2 &1.3E+05\\ 
     1 &   1E+03 &  1E+03 &     0 & 1.8E-03 & 9.7E-04 & 2.1E-02 & 1.3E-02 &      2 &1.3E+05\\ 
     1 &   1E+03 &  1E+04 &     0 & 1.8E-03 & 9.7E-04 & 2.1E-02 & 1.3E-02 &      2 &1.3E+05\\ 
     1 &   1E+04 &  1E+02 &     0 & 1.6E-10 & 8.3E-06 & 7.1E-06 & 1.3E-03 &      1 &6.5E+04\\ 
     1 &   1E+04 &  1E+03 &     0 & 1.6E-10 & 8.3E-06 & 7.1E-06 & 1.3E-03 &      1 &6.5E+04\\ 
     1 &   1E+04 &  1E+04 &     0 & 1.6E-10 & 8.3E-06 & 7.1E-06 & 1.3E-03 &      1 &6.5E+04\\ 
     1 &   1E+04 &  1E+05 &     0 & 1.6E-10 & 8.3E-06 & 7.1E-06 & 1.3E-03 &      1 &6.5E+04\\ 
     \hline 
     2 &   1E+02 &  1E+02 &     0 & 1.7E-01 & 1.0E-01 & 2.0E-01 & 1.2E-01 &      2 &1.3E+05\\ 
     2 &   1E+02 &  1E+03 &     0 & 1.7E-01 & 1.0E-01 & 2.0E-01 & 1.2E-01 &      2 &1.3E+05\\ 
     2 &   1E+03 &  1E+02 &     0 & 1.8E-03 & 9.7E-04 & 2.1E-02 & 1.3E-02 &      2 &1.3E+05\\ 
     2 &   1E+03 &  1E+03 &     0 & 1.8E-03 & 9.7E-04 & 2.1E-02 & 1.3E-02 &      2 &1.3E+05\\ 
     2 &   1E+03 &  1E+04 &     0 & 1.8E-03 & 9.7E-04 & 2.1E-02 & 1.3E-02 &      2 &1.3E+05\\ 
     2 &   1E+04 &  1E+02 &     0 & 1.6E-10 & 8.3E-06 & 7.1E-06 & 1.3E-03 &      1 &6.5E+04\\ 
     2 &   1E+04 &  1E+03 &     0 & 1.6E-10 & 8.3E-06 & 7.1E-06 & 1.3E-03 &      1 &6.5E+04\\ 
     2 &   1E+04 &  1E+04 &     0 & 1.6E-10 & 8.3E-06 & 7.1E-06 & 1.3E-03 &      1 &6.5E+04\\ 
     2 &   1E+04 &  1E+05 &     0 & 1.6E-10 & 8.3E-06 & 7.1E-06 & 1.3E-03 &      1 &6.5E+04\\ 
     \hline 
     3 &  1E+02 & -     &     0 & 2.5E-03 & 1.3E-03 & 2.4E-02 & 1.5E-02 &      2 &1.3E+05\\ 
     3 &  1E+03 & -     &     0 & 6.4E-10 & 1.3E-05 & 1.4E-05 & 1.6E-03 &      1 &6.6E+04\\ 
     3 &  1E+04 & -     &     0 & 6.3E-14 & 1.3E-07 & 1.5E-07 & 1.6E-04 &      1 &6.7E+04\\ 
    \end{tabular}
    \caption{Parallel free surface to the reference slope test problem for $n_{el}=100$ and $\sigmao=\SI{9.035}{Pa}$.}
    \label{tab:errorsEx2}
\end{table}

For all the cases in this study, the fluid is 0\% active, i.e., rigid in the entire domain at time $t=1$.
On average, it takes 1 or 2 Newton-Raphson iterations for each time step (see column NRIt/Tit in Table \ref{tab:errorsEx2} ) for all the regularization approaches.
This convergence is quadratic which is characteristic of the Newton-Raphson method.
Note that for all the approaches, the error for the height is smaller than the error for the velocity.
The height and the zero velocity are better preserved, i.e., lower errors, for larger values of the parameters $\gamma$ and $\beta$.
Note that the third approach is more accurate for this example using the same parameter $\beta$, and the computational time is smaller than the other two approaches.
We attribute this to the smoothness of the third regularization function, which has a continuous first derivative.

\FloatBarrier
\subsection{Dam break}

In this example, we model a reservoir with  
fluid level $h_1=\SI{1.5}{\meter}$ higher than the horizontal ground. 
A dam wall at $x=\SI{1.5}{\meter}$ is removed instantaneously to simulate the dam break.
The initial velocity and auxiliary variable $\EE$ are both set to zero.
In this case, the dam break flows into a wet channel of height $h_2=\SI{0.5}{\meter}$.
We consider a domain of length $\Omega=[0,L]$ where $L=3$ and apply boundary conditions on both sides as follows $h(0,t)=h_1$, $h(L,t)=h_2$, $\partial\vel(0,t)/\partial x=\EE(0,t)=\partial\vel(L,t)/\partial x=\EE(L,t)=0$.

\begin{figure}
    \centering
    \includegraphics{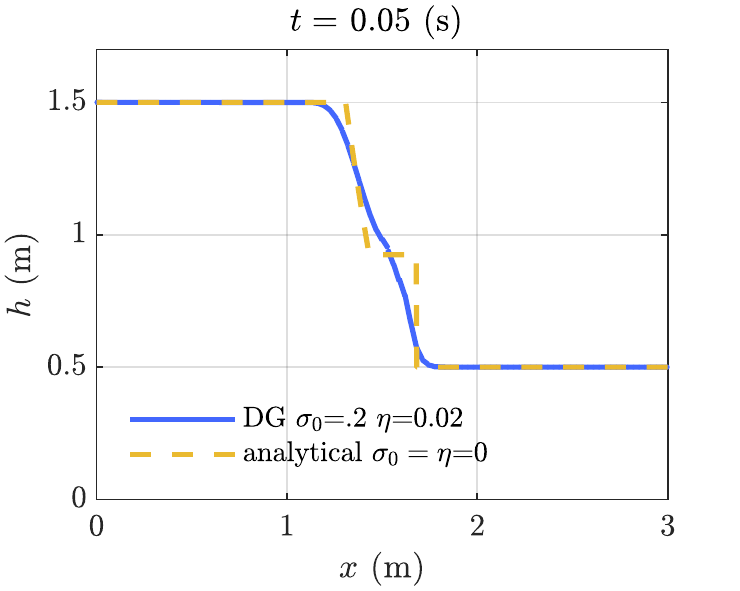}
    \includegraphics{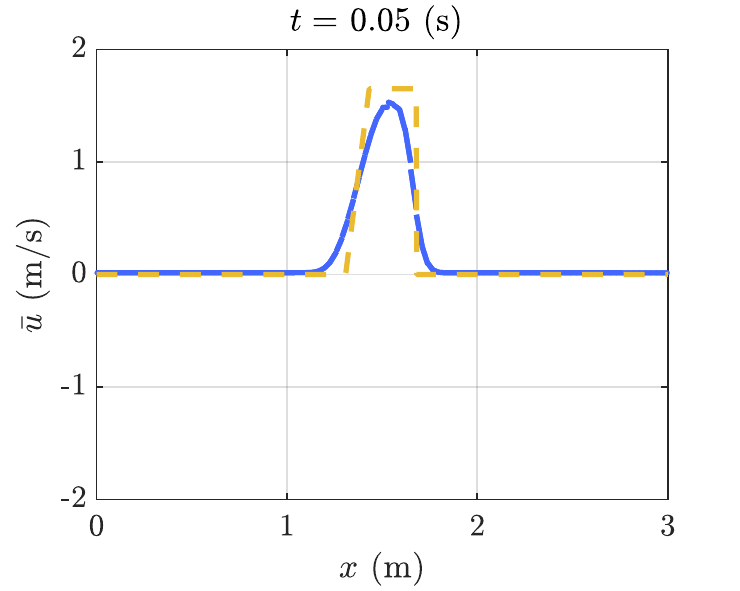}
    \includegraphics{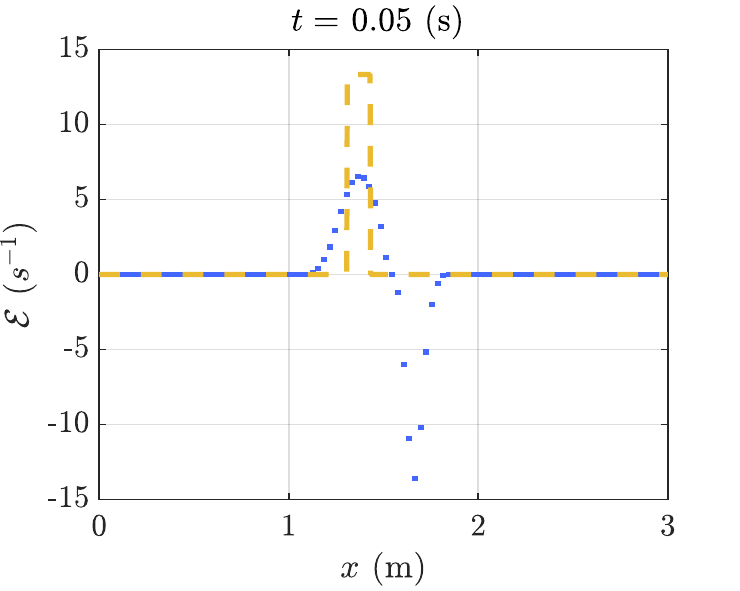}
     \includegraphics{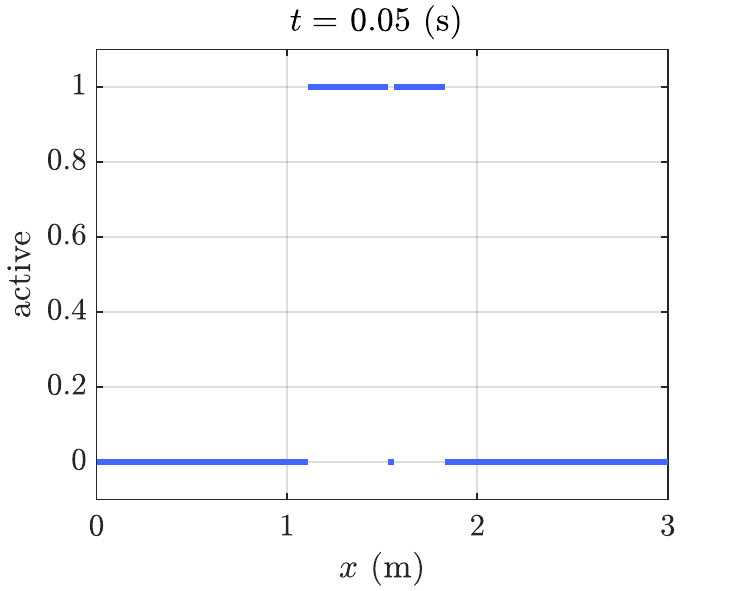}
    \caption{Discontinous Galerkin 100 elements including viscosity and Bingham constitutive law for the dam break problem on a wet bed at time $t=\SI{0.05}{\second}$ with $\sigmao=\SI{0.2}{\pascal}$ and $\eta=\SI{0.02}{\pascal\second}$.}
    \label{Fig:Dambreak1}
\end{figure}

\begin{figure}
    \centering
    \includegraphics{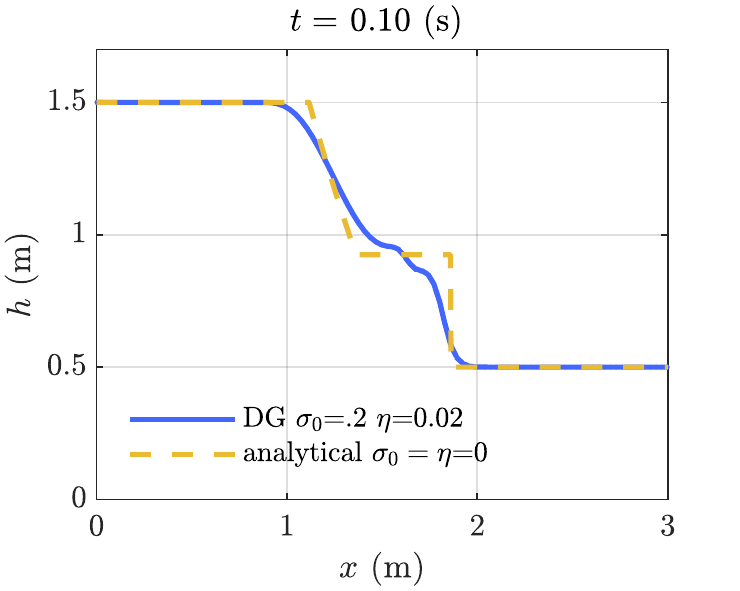}
    \includegraphics{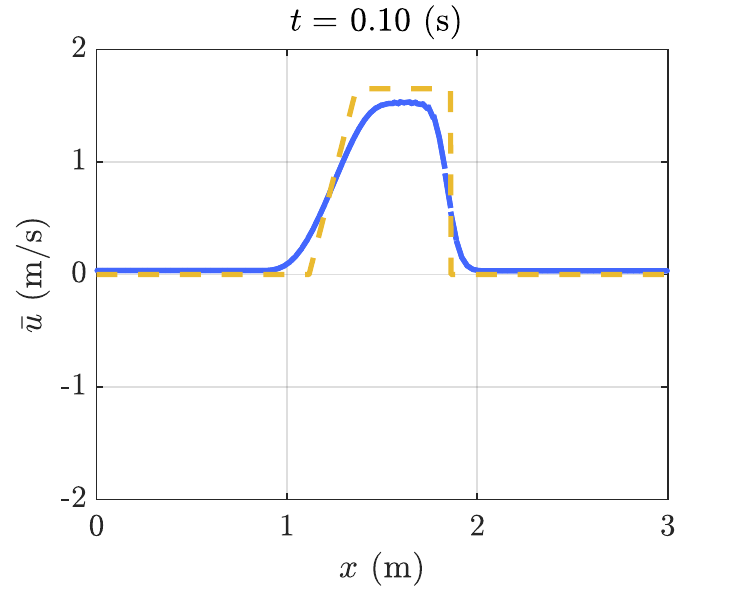}
    \includegraphics{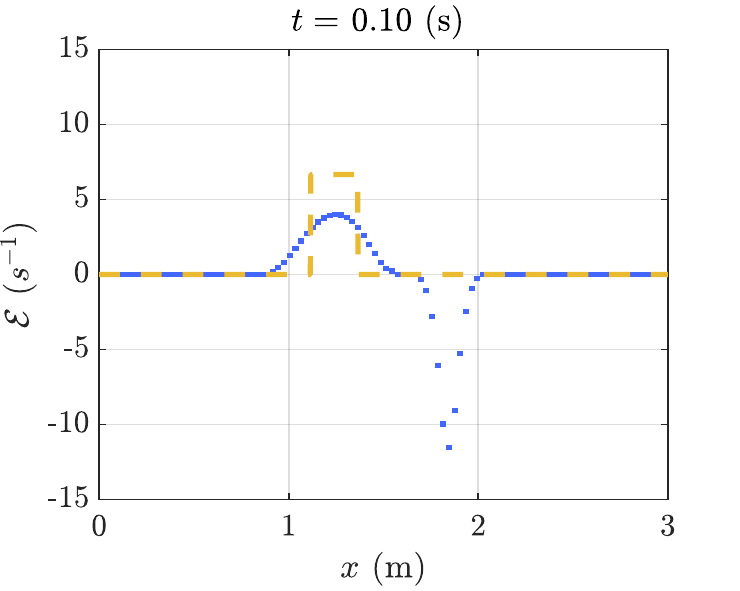}
     \includegraphics{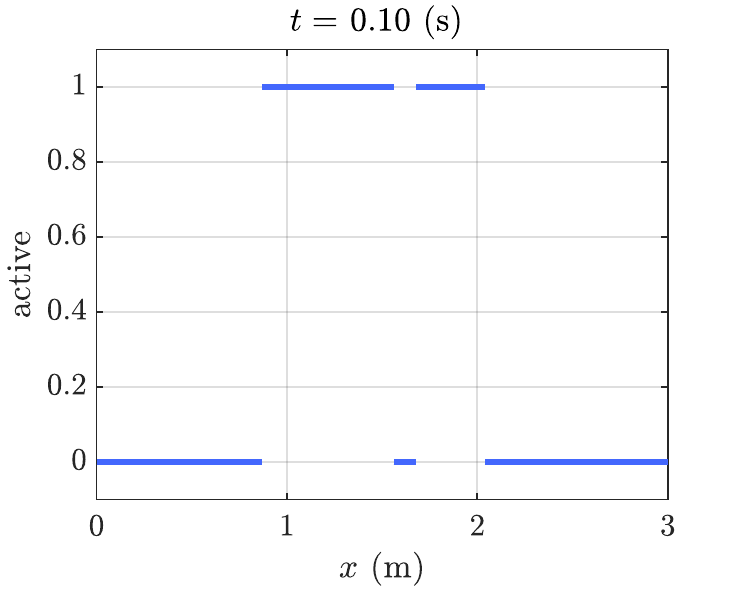}
    \caption{Discontinous Galerkin 100 elements including viscosity and Bingham constitutive law for the dam break problem on a wet bed at time $t=\SI{0.10}{\second}$ with $\sigmao=\SI{0.2}{\pascal}$ and $\eta=\SI{0.02}{\pascal\second}$.}
    \label{Fig:Dambreak2}
\end{figure}

\begin{figure}
    \centering
    \includegraphics{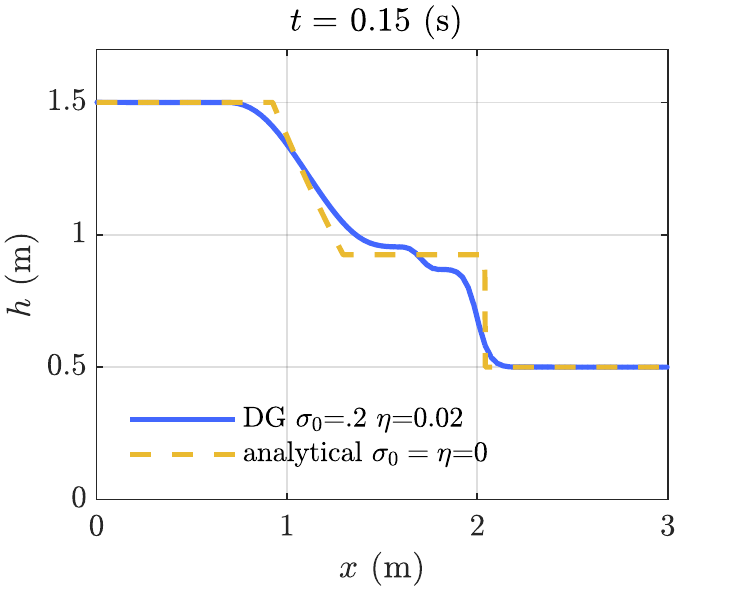}
    \includegraphics{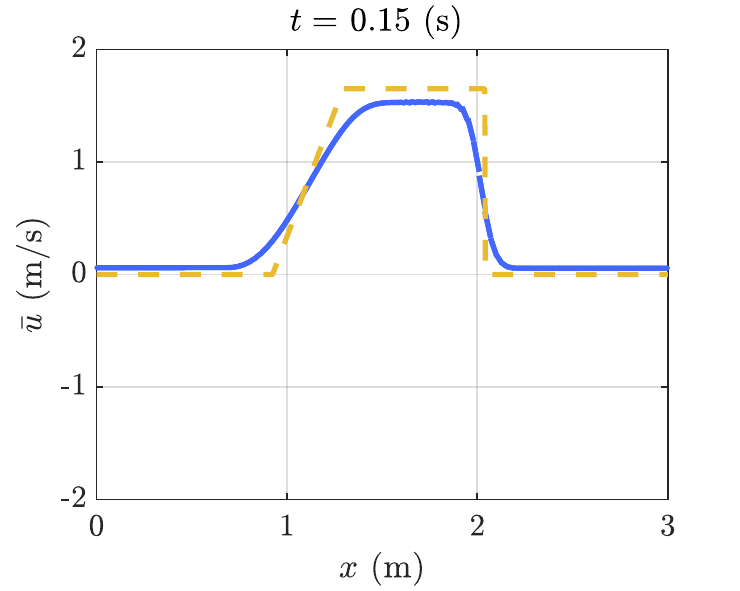}
    \includegraphics{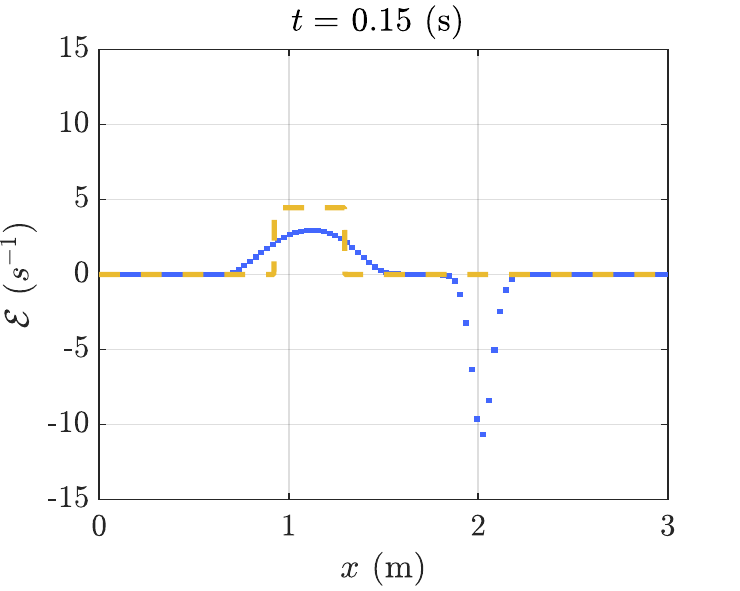}
     \includegraphics{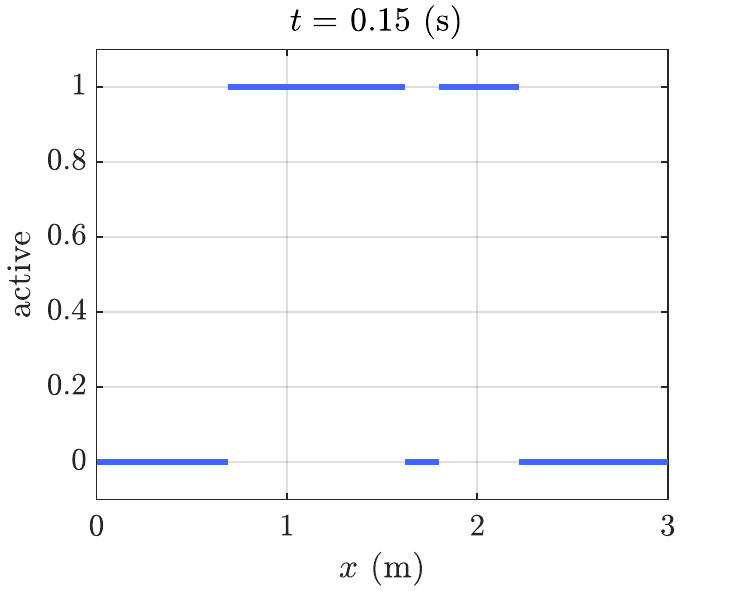}
    \caption{Discontinous Galerkin 100 elements including viscosity and Bingham constitutive law for the dam break problem on a wet bed at time $t=\SI{0.15}{\second}$ with $\sigmao=\SI{0.2}{\pascal}$ and $\eta=\SI{0.02}{\pascal\second}$.}
    \label{Fig:Dambreak3}
\end{figure}

We solve this dam-break problem using 100 elements, linear basis functions for the height $h$ and the velocity $\vel$, and piece wise constants for the auxiliary variable $\EE$, i.e.,  $m_1=m_2=1$ and $m_3=0$.
The Bingham yield stress, regularization parameters, gravity, density, and viscosity are  $\sigmao=\SI{0.2}{\pascal}$, $\gamma = 10^2\SI{}{\pascal\second}$, $\beta = 10^2\SI{}{\pascal\second}$, $g=\SI{9.81}{\meter/\second^2}$, $\rho=\SI{1}{kg.m^{-3}}$, and $\eta=\SI{0.02}{\pascal\second}$, respectively.
The time step for the computation is $\Delta t = 10^{-6}\SI{}{s}$.
We show the solution for this problem in Figures \ref{Fig:Dambreak1}, \ref{Fig:Dambreak2}, and \ref{Fig:Dambreak3} at times = 0.05, 0.10, 0.15 s respectively.
To compute the Bingham stress contribution, we use the first alternative regularization approach with piece wise continuous derivative, cf. Equation \eqref{Eq:sigmaTwoA}.
If we use the other alternatives, i.e., \eqref{Eq:sigmaTwoB} and \eqref{Eq:sigmaTwoC}, the figures are very similar.
The variable $\EE$ representing the gradient of the velocity, i.e., strain rate, is shown in Figures \ref{Fig:Dambreak1}c, \ref{Fig:Dambreak2}c, and \ref{Fig:Dambreak3}c.
With this variable $\EE$, it is easy to determine the active  regions such $|\EE|>\sigmao/\gamma$, cf. Figures \ref{Fig:Dambreak1}d, \ref{Fig:Dambreak2}d, and \ref{Fig:Dambreak3}d.
For reference and comparison, the analytical solution of this problem with no viscosity and no Bingham yield stress, i.e., $\eta=0$ and $\sigmao=0$, is presented in the figures, which solution is detailed in \cite{Stoker2011}.

\begin{figure}
    \centering
    \includegraphics{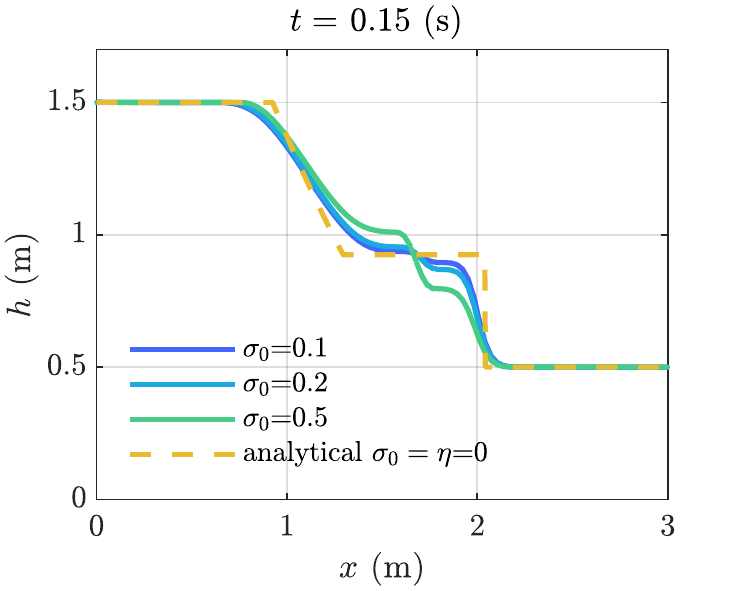}
    \includegraphics{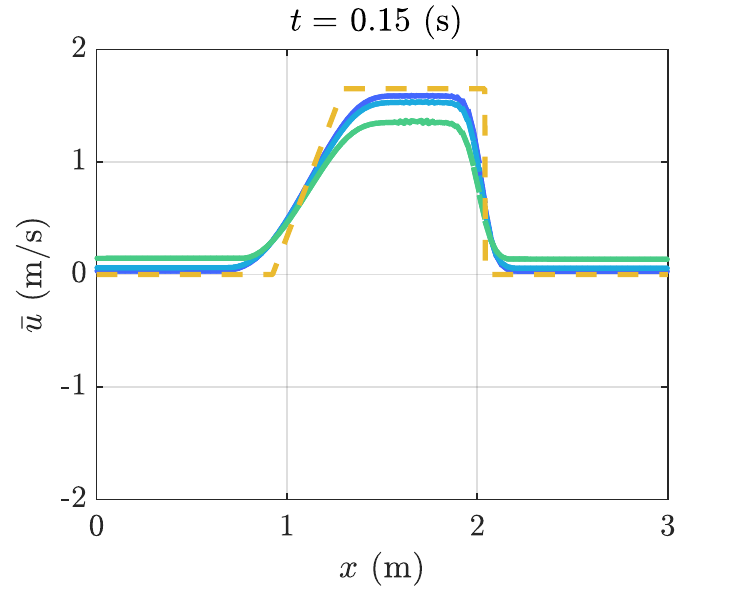}
    \caption{Discontinous Galerkin 100 elements including viscosity and Bingham constitutive law for the dam break problem on a wet bed at time $t=$0.15 s for different Bingham yield stress $\sigma_0$ with viscosity $\eta=.02$.}
    \label{Fig:Dambreak3sigma}
\end{figure}

\begin{figure}
    \centering
    \includegraphics{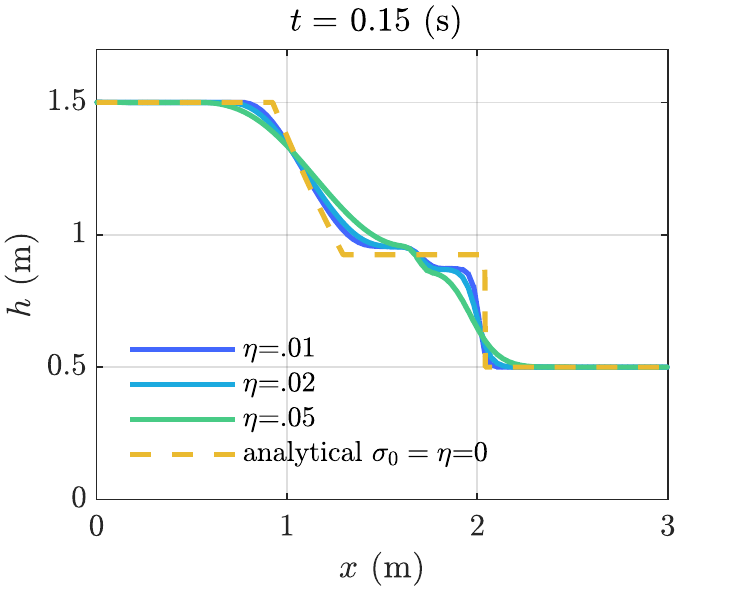}
    \includegraphics{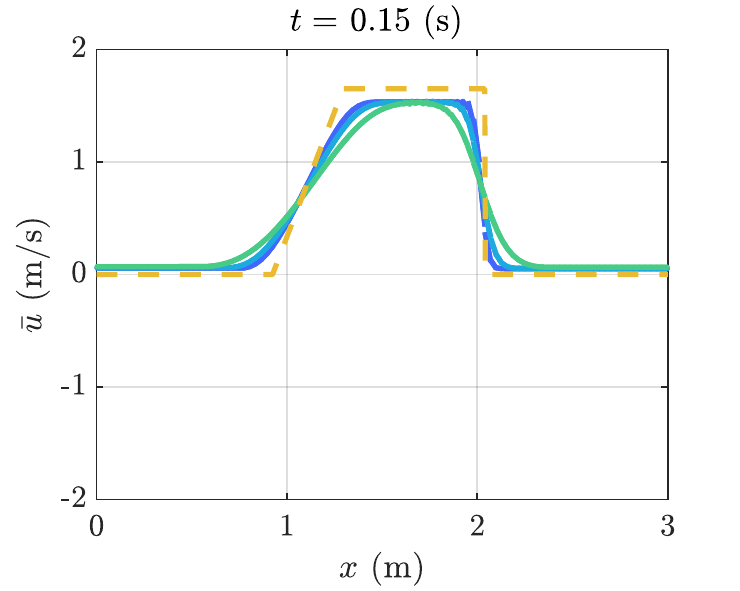}
    \caption{Discontinous Galerkin 100 elements including viscosity and Bingham constitutive law for the dam break problem on a wet bed at time $t=$0.15 s with Bingham yield stress $\sigma_0=.2$ for different viscosity $\eta$.}
    \label{Fig:Dambreak3nu}
\end{figure}

Due to the viscosity of the material, the solution is smooth despite the discontinuity in the initial condition as we can see this effect in Figure \ref{Fig:Dambreak3nu}.
Indeed, the flow resistance of the viscosity smooths spurious oscillations, and we do not need extra computational techniques, e.g., slope limiters, to deal with these oscillations.
Due to the Bingham viscoplasticity, the domain has five regions, i.e., inactive-left, active-left-centered, inactive-centered, active-right-centered, and inactive-right.
The central region consists of the active-left-centered, inactive-centered, and active-right-centered portions, and it grows as time passes.
An important effect is that if the Bingham yield stress is large, the velocity for the central region gets lower as shown in Figure \ref{Fig:Dambreak3sigma}.
In the inactive-centered region, we observe another step in height because the material behaves like a solid there.
This centered step is not present if the Bingham yield stress is zero, and it gets bigger for higher Bingham yield stresses.
We also note that the inactive-left and inactive-right regions acquire nonzero velocities but move like solids.

\FloatBarrier

We repeat the above numerical experiment with $n_{ele}=100$ elements, $\sigmao=\SI{0.2}{\pascal}$ and $\eta=\SI{0.02}{\pascal\second}$ for the different regularization approaches.
Since we do not have an analytical solution, we compare our results with a numerical experiment with a finer mesh of $n_{ele}=500$ and smaller time steps $\Delta t = 10^{-6}\SI{}{s}$, cf. Table \ref{Tab:refined}.
We compute the error norms L-2 and L-infinity comparing these solutions to the refined one.
In Table \ref{tab:damBreakFuns}, we summarize the results.
First, we note that for regularization approaches 1 and 2 with large values of $\gamma=1e3$, we have convergence difficulties, high computational costs, and low accuracy.
In contrast, for regularization parameters $\gamma<1e3$, the regularization approach 1 and 2 present low computational cost and great accuracy.
Additionally, the prediction of the percentage of active regions is close to the 46\% active that the refined solution predicts, so approaches 1 and 2 have the closest values compared to the other methods.
Also, we do not find any significant effect on the regularization parameter $\beta$.
We obtain a quadratic convergence at each time step as expected for the Newton-Raphson method when we use the regularization approach 1 and 2 with $\gamma<1e3$ and approach 3 for any regularization parameter.
This clearly shows that the smoothness of the regularization function affect convergence.
An important detail is that we can overcome these convergence difficulties using smaller time steps which is the case of our refined solution where we do not have convergence issues.
Both approaches 1 and 2 converge with 2 Newton-Raphson iterations for the refined case, cf. Table \ref{Tab:refined}, and again obtain a quadratic convergence.
About the third regularization approach, the computational cost is low and accuracy is great for any regularization parameter.
We attribute this cost-effectiveness to the $C^\infty$ smoothness that helps the Newton-Raphson convergence.
However, this approach 3 fails to predict the active percent region because the regularization function does not possess a clear division of active-inactive region due to its natural smoothness.
Despite the error to predict the active regions, approach 3 has a low error in height and velocity compared to the refined solution.

\begin{table}[!ht]
    \centering
    \begin{tabular}{cccccc}
         Reg. &  $\gamma$ &  $\beta$ &  active &  NRIt/TIt & CPUTime \\
          approach & $(\SI{}{\pascal\second})$  &   & $\%$   & & (s) \\
    \hline
     1 &   1E+03 &  1E+03 &   46.0 &    2.0 &6.8E+05 \\ 
     2 &   1E+03 &  1E+03 &   46.0 &    2.0 &6.8E+05 \\ 
     3 &  1E+03 & -     &   90.2 &    2.0 &6.8E+05  \\  
 \hline 
    \end{tabular}
    \caption{Refined solution for dam break problem with $n_{el}=500$, $\Delta t = 10^{-6}\SI{}{s}$, $\sigmao=0.2\SI{}{Pa}$, and $\eta=0.02\SI{}{\pascal\second}$ using different regularization functions.}
    \label{Tab:refined}
\end{table}

\begin{table}[!ht]
    \centering
    \begin{tabular}{cccccccccc}
         Reg. &  $\gamma$ &  $\beta$ &  active &  \Lth & \Ltvel & \Lh & \Lvel & NRIt/TIt & CPUTime \\
          approach & $(\SI{}{\pascal\second})$  &   & $\%$ &   &  &  &  &  & (s) \\
    \hline
     1 &   1E+01 &  1E+01 &   48.0 &2.2E-04 &2.0E-03 &3.1E-02 &1.3E-01 &    2.0 &1.4E+04 \\ 
     1 &   1E+01 &  1E+02 &   44.0 &2.2E-04 &2.0E-03 &3.1E-02 &1.3E-01 &    2.0 &1.4E+04 \\ 
     1 &   1E+02 &  1E+02 &   45.0 &2.2E-04 &2.2E-03 &3.2E-02 &1.4E-01 &    2.0 &1.4E+04 \\ 
     1 &   1E+02 &  1E+03 &   45.0 &2.2E-04 &2.2E-03 &3.2E-02 &1.4E-01 &    2.0 &1.4E+04 \\ 
     1 &   1E+03 &  1E+02 &   12.0 &1.8E-01 &2.1E+00 &5.4E-01 &1.5E+00 &   10.0 &7.0E+04 \\ 
     1 &   1E+03 &  1E+03 &   12.0 &1.8E-01 &2.1E+00 &5.4E-01 &1.5E+00 &   10.0 &7.0E+04 \\ 
     1 &   1E+03 &  1E+04 &   12.0 &1.8E-01 &2.1E+00 &5.4E-01 &1.5E+00 &   10.0 &7.0E+04 \\ 
     \hline 
     2 &   1E+01 &  1E+01 &   48.0 &2.2E-04 &2.0E-03 &3.1E-02 &1.3E-01 &    2.0 &1.4E+04 \\ 
     2 &   1E+01 &  1E+02 &   44.0 &2.2E-04 &2.0E-03 &3.1E-02 &1.3E-01 &    2.0 &1.4E+04 \\ 
     2 &   1E+02 &  1E+02 &   45.0 &2.2E-04 &2.2E-03 &3.2E-02 &1.4E-01 &    2.0 &1.4E+04 \\ 
     2 &   1E+02 &  1E+03 &   45.0 &2.2E-04 &2.2E-03 &3.2E-02 &1.4E-01 &    2.0 &1.4E+04 \\ 
     2 &   1E+03 &  1E+02 &   12.0 &1.8E-01 &2.1E+00 &5.4E-01 &1.5E+00 &   10.0 &7.0E+04 \\ 
     2 &   1E+03 &  1E+03 &   12.0 &1.8E-01 &2.1E+00 &5.4E-01 &1.5E+00 &   10.0 &7.0E+04 \\ 
     2 &   1E+03 &  1E+04 &   12.0 &1.8E-01 &2.1E+00 &5.4E-01 &1.5E+00 &   10.0 &7.0E+04 \\ 
     \hline 
     3 &  1E+01 & -     &   88.0 &2.6E-04 &2.3E-03 &3.2E-02 &1.1E-01 &    2.0 &1.4E+04  \\ 
     3 &  1E+02 & -     &   90.0 &2.2E-04 &2.1E-03 &3.1E-02 &1.3E-01 &    2.0 &1.4E+04  \\ 
     3 &  1E+03 & -     &   90.0 &2.2E-04 &2.2E-03 &3.3E-02 &1.4E-01 &    2.0 &1.4E+04  \\ 
     3 &  1E+04 & -     &   90.0 &2.2E-04 &2.2E-03 &3.3E-02 &1.4E-01 &    2.0 &1.4E+04  \\  
     \hline
    \end{tabular}
    \caption{Dam break problem for $n_{el}=100$, $\sigmao=0.2\SI{}{Pa}$, $\eta=0.02\SI{}{\pascal\second}$, $m_1=m_2=1$, $m_3=0$, and using different regularization functions and parameters.}
    \label{tab:damBreakFuns}
\end{table}

We repeat the same example as before, but we change the polynomial order of the auxiliary variable, i.e., $m_3=1$.
Table \ref{tab:damBreakFuns1} shows the summary of these results.
We have the same behavior as before with all the methods.
There is a slightly better improvement in accuracy for approaches 1 and 2 with $\gamma<1e3$.
Similarly, approach 3 has better accuracy for $\gamma<1e4$, and for $\gamma=1e4$, there are convergence issues and the accuracy error is larger.
Finally, for all the methods, the computational cost is higher with $m_3=1$ than with $m_3=0$ as expected.

\begin{table}[!ht]
    \centering
    \begin{tabular}{cccccccccc}
         Reg. &  $\gamma$ &  $\beta$ &  active &  \Lth & \Ltvel & \Lh & \Lvel & NRIt/TIt & CPUTime \\
          approach & $(\SI{}{\pascal\second})$  &   & $\%$ &   &  &  &  &  & (s) \\
    \hline
     1 &   1E+01 &  1E+01 &   49.8 &2.0E-04 &1.9E-03 &3.0E-02 &1.2E-01 &    2.0 &1.9E+04 \\ 
     1 &   1E+02 &  1E+02 &   48.7 &1.8E-04 &2.0E-03 &3.0E-02 &1.2E-01 &    2.1 &2.0E+04 \\ 
     1 &   1E+03 &  1E+03 &   53.1 &1.7E-01 &2.0E+00 &5.3E-01 &1.5E+00 &   10.0 &9.5E+04 \\ 
     \hline 
     2 &   1E+01 &  1E+01 &   50.1 &2.0E-04 &1.9E-03 &3.0E-02 &1.2E-01 &    2.0 &1.9E+04 \\ 
     2 &   1E+02 &  1E+02 &   48.7 &1.8E-04 &2.0E-03 &3.0E-02 &1.2E-01 &    2.1 &2.0E+04 \\ 
     2 &   1E+03 &  1E+03 &   55.6 &1.7E-01 &2.0E+00 &5.3E-01 &1.5E+00 &   10.0 &9.5E+04 \\ 
     \hline 
     3 &  1E+01 & -     &   88.5 &2.5E-04 &2.3E-03 &2.8E-02 &1.0E-01 &    2.0 &1.9E+04  \\ 
     3 &  1E+02 & -     &   90.8 &2.0E-04 &1.9E-03 &3.0E-02 &1.2E-01 &    2.1 &1.9E+04  \\ 
     3 &  1E+03 & -     &   89.4 &1.9E-04 &2.0E-03 &3.0E-02 &1.2E-01 &    2.1 &2.0E+04  \\ 
     3 &  1E+04 & -     &   82.3 &1.7E-01 &2.0E+00 &5.3E-01 &1.5E+00 &   10.0 &9.5E+04  \\ 
     \hline 
 \end{tabular}
    \caption{Dam break problem for $n_{el}=100$, $\sigmao=0.2\SI{}{Pa}$, $\eta=0.02\SI{}{\pascal\second}$, , $m_1=m_2=m_3=1$, and using different regularization functions and parameters.}
    \label{tab:damBreakFuns1}
\end{table}

To reduce the convergence issues of approaches 1, and 2, we can use smaller time steps, or another alternative is the \emph{continuation} method.
In each time step, we can solve the problem for a regularization parameter $\gamma_i$, achieve convergence, then increase the value of the regularization parameter to $\gamma_{i+1}$.
We iterate this procedure until we reach the desired $\gamma$.
In our numerical experiment, we start from a regularization parameter $\gamma_0$ and reach $\gamma$ in $n_\gamma\ge2$ continuation steps.
So, in each continuation iteration $i=1,2,..., n_\gamma$ the regularization parameter is $\gamma_i= \gamma_0 + (i-1)(\gamma -\gamma_0)/(n_\gamma-1)$.
To test this method, we repeat the example for $n_{el}=100$, $\sigmao=0.2\SI{}{Pa}$, $\eta=0.02\SI{}{\pascal\second}$, $m_1=m_2=1$, $m_3=0$, and using the regularization approach 1 and 2.
We present these results in Table \ref{tab:damBreakContinuation}.
For this test, we start and finish with $\gamma_0=10^2$ and $\gamma=10^3$.
In all these cases, the continuation method successfully converged with great accuracy.
As expected, the computational cost is higher than the results with no continuation.
Interestingly, we obtained the lowest computational cost for $n_\gamma=2$ continuation iterations.
As we can see, the continuation method for the regularization parameter is a good alternative to avoid convergence issues.

\begin{table}[!ht]
    \centering
    \begin{tabular}{ccccccccccc}
         Reg. &  $\gamma_0$ & $\gamma$ &  $n_\gamma$ &  active &  \Lth & \Ltvel & \Lh & \Lvel & NRIt/TIt & CPUTime \\
          approach &$(\SI{}{\pascal\second})$ & $(\SI{}{\pascal\second})$  &   & $\%$ &   &  &  &  &  & (s) \\
    \hline
     1 &   1E+02 &  1E+03 &   2 &   94.0 &2.7E-03 &3.3E-02 &1.3E-01 &6.4E-01 &   11.5 &8.0E+04 \\ 
     1 &   1E+02 &  1E+03 &   3 &   28.0 &2.7E-03 &3.3E-02 &1.3E-01 &6.4E-01 &   11.8 &8.2E+04 \\ 
     1 &   1E+02 &  1E+03 &   5 &   21.0 &2.7E-03 &3.3E-02 &1.3E-01 &6.4E-01 &   14.5 &1.0E+05 \\ 
     1 &   1E+02 &  1E+03 &  10 &   21.0 &2.7E-03 &3.3E-02 &1.3E-01 &6.4E-01 &   22.3 &1.5E+05 \\ 
     \hline 
     2 &   1E+02 &  1E+03 &   2 &   94.0 &2.7E-03 &3.3E-02 &1.3E-01 &6.4E-01 &   11.5 &8.0E+04 \\ 
     2 &   1E+02 &  1E+03 &   3 &   28.0 &2.7E-03 &3.3E-02 &1.3E-01 &6.4E-01 &   11.8 &8.3E+04 \\ 
     2 &   1E+02 &  1E+03 &   5 &   21.0 &2.7E-03 &3.3E-02 &1.3E-01 &6.4E-01 &   14.5 &1.0E+05 \\ 
     2 &   1E+02 &  1E+03 &  10 &   21.0 &2.7E-03 &3.3E-02 &1.3E-01 &6.4E-01 &   22.3 &1.6E+05 \\ 
     \hline 
\end{tabular}
    \caption{Dam break problem for $n_{el}=100$, $\sigmao=0.2\SI{}{Pa}$, $\eta=0.02\SI{}{\pascal\second}$, and using a continuation method for the regularization parameter.}
    \label{tab:damBreakContinuation}
\end{table}

We repeat the example for $n_{el}=50$, $\sigmao=0.2\SI{}{Pa}$, $\eta=0.02\SI{}{\pascal\second}$, $m_1=m_2=1$, $m_3=0$, and use in this case different polynomial orders for the basis of the height, velocity, and auxiliary variable.
We obtain good results for all the approaches.
If we use a higher-order polynomial for the basis, we get lower errors but higher computational times, as expected.
Approaches 1 and 2 predict the active percent region with better accuracy than approach 3.

\begin{table}[!ht]
    \centering
    \begin{tabular}{ccccccccccc}
         Reg. &  $m_1$ & $m_2$ &  $m_3$ &  active &  \Lth & \Ltvel & \Lh & \Lvel & NRIt/TIt & CPUTime \\
          approach &$(\SI{}{\pascal\second})$ & $(\SI{}{\pascal\second})$  &   & $\%$ &   &  &  &  &  & (s) \\
    \hline
     1 &      2 &      2 &       2 &   49.2 &3.4E-04 &3.5E-03 &8.4E-02 &3.7E-01 &    2.0 &1.9E+05 \\ 
     1 &      3 &      3 &       3 &   50.2 &2.2E-04 &2.3E-03 &7.6E-02 &3.1E-01 &    2.0 &3.1E+05 \\ 
     1 &      4 &      4 &       4 &   50.9 &2.0E-04 &1.9E-03 &7.2E-02 &3.0E-01 &    2.0 &4.9E+05 \\ 
     \hline 
     2 &      2 &      2 &      2  &   49.2 &3.4E-04 &3.5E-03 &8.4E-02 &3.7E-01 &    2.0 &1.9E+05 \\ 
     2 &      3 &      3 &      3  &   50.2 &2.2E-04 &2.3E-03 &7.6E-02 &3.1E-01 &    2.0 &3.1E+05 \\ 
     2 &      4 &      4 &      4  &   53.2 &2.1E-04 &1.9E-03 &7.2E-02 &3.0E-01 &    2.0 &5.0E+05 \\ 
     \hline 
     3 &      2 &      2 &      2 &   91.9 &3.5E-04 &3.3E-03 &8.4E-02 &3.7E-01 &    2.0 &1.9E+05  \\ 
     3 &      3 &      3 &      3 &   86.8 &2.3E-04 &2.2E-03 &7.5E-02 &3.1E-01 &    2.0 &3.1E+05  \\ 
     3 &      4 &      4 &      4 &   90.3 &2.0E-04 &1.8E-03 &7.1E-02 &3.0E-01 &    2.0 &4.9E+05  \\ 
     \hline 
     \end{tabular}
    \caption{Dam break problem for $n_{el}=50$, $\sigmao=0.2\SI{}{Pa}$, $\eta=0.02\SI{}{\pascal\second}$, and using different polynomial orders for basis.}
    \label{tab:damBreakmodes}
\end{table}

\FloatBarrier
\section{Conclusions}
\label{Sec:Conclusions}
This paper simulates viscoplastic flow by solving the shallow-water equations using the discontinuous Galerkin method.
Due to the element-base discretization, the discontinuous Galerkin method is attractive due to its high parallelization, h- and p-adaptivity, and ability to capture the discontinuities of the exact solution.
Since cell interfaces are shared by two elements, we use numerical fluxes based on the theory of Riemann problems to ensure a stable solution of the nonlinear hyperbolic equations. 
For the viscoplastic material model, we use a depth-averaged Bingham constitutive law in which the material behaves as a solid or fluid depending if the stress magnitude is below or above the material's yield stress.
To couple the Bingham model with the shallow-water equations, we regularize this constitutive law~\cite{DeLosReyes2012a}.
The idea behind the regularization is that the material presents a large viscosity if the stress is below the yield threshold; and in this way, the behavior mimics a solid.
On the other hand, the material presents the usual viscosity if the stress is above the yield threshold.
We study three regularization functions, two of local type and one global, with continuous first derivatives.
These regularizations yield a system of Newton's differentiable equations where we obtain quadratic convergence in each time iteration solution.
Numerical examples tested well-balanced properties and the effectiveness of our approach.
We compared the different regularization approaches proposed.
And we found that the $C^\infty$ regularization is computationally more efficient, however, the method cannot predict accurately the active/inactive regions.
In contrast, the other piece-wise linear (local) regularizations predict better the active/inactive regions.
To avoid convergence issues, we propose to use a continuation method on the regularization parameter or smaller time steps.

In the current work. we solved one-dimensional problems in 1D. So for future work, we will extend our approach to two-dimensional problems and include treatment for dry/wet fronts.

\section{Acknowledgments}
This research was carried out using the research computing facilities and/or advisory services offered by the Scientific Computing Laboratory of the Research Center on Mathematical Modeling: MODEMAT, Escuela Politécnica Nacional - Quito.
This project was funded by the Escuela Politécnica Nacional through the project PIM 20-01.

\section*{Appendix}
\setcounter{section}{0}
\setcounter{equation}{0}
\setcounter{figure}{0}
\renewcommand{\theequation}{\thesection.\arabic{equation}}
\renewcommand{\thefigure}{\thesection.\arabic{figure}}
\renewcommand{\thesection}{\Alph{section}}
\setcounter{equation}{0}
\setcounter{figure}{0}
\section{Tangent terms}
\label{Ape:Derivatives}
In this section, we detail the derivatives for the residual vectors of Equation \eqref{Eq:WeakForm7}.
The derivatives are essential to building the tangent terms for our Newton System \eqref{Eq:Newton}. 
We systematically follow the chain rule and obtain the following expressions for the element $\Omega_e^h$:
\begin{subequations}
    \begin{align}
        \derivev{ \bfR_{Ve}^{(k)} }&=  
        \frac{1}{\Delta t}\Bigg\langle \bfNV 
        \, ,  \derivev{ \bfU^{(k)} }  \Bigg\rangle_{\Omega_e^h} \,
          + \Bigg\langle \bfNV^- \,,  \derivev{\nbfF^{(k)}} \, \cdot \bfne \Bigg\rangle_{\Gamma_e^h} 
         - \Bigg\langle \derive{ \bfNV}{x_1} \,, \derivev{ \bfF_1^{(k)} } \Bigg\rangle_{\Omega_e^h}
          - \Bigg\langle \derive{ \bfNV}{x_2} \,, \derivev{ \bfF_2^{(k)} } \Bigg\rangle_{\Omega_e^h}
          \nonumber \\& \quad 
         -\Bigg\langle \bfNV^-  \,,  \derivev{ \nbfQ^{(k)} } \, \cdot \bfne \Bigg\rangle_{\Gamma_e^h}
        +\Bigg\langle \derive{\bfNV}{x_1} \,, \derivev{ \bfQ_1^{(k)} }\Bigg\rangle_{\Omega_e^h}
        +\Bigg\langle \derive{\bfNV}{x_2} \,, \derivev{ \bfQ_2^{(k)} }\Bigg\rangle_{\Omega_e^h}
         \nonumber \\& \quad 
         -\Bigg\langle \bfNV \,, \derivev{ \bfS^{(k)}} \Bigg\rangle_{\Omega_e^h} 
         \\
         \derivee{ \bfR_{Ve}^{(k)} }&=  
         -\Bigg\langle \bfNV^-  \,,  \derivee{ \nbfQ^{(k)} } \, \cdot \bfne \Bigg\rangle_{\Gamma_e^h}
        +\Bigg\langle \derive{ \bfNV}{x_1} \,, \derivee{ \bfQ_1^{(k)} }\Bigg\rangle_{\Omega_e^h}
        +\Bigg\langle \derive{ \bfNV}{x_2} \,, \derivee{ \bfQ_2^{(k)} }\Bigg\rangle_{\Omega_e^h}
        \\
        \derivev{ \bfR_{Ee}^{(k)} }&= 
         -\Bigg\langle \bfNE^-  \,,  \derivev{ \nbfG^{(k)}} \, \cdot \bfne \Bigg\rangle_{\Gamma_e^h}
        +\Bigg\langle \frac{ \partial \bfNE}{\partial x_1} \,, \derivev{ \bfG_1^{(k)} } \Bigg\rangle_{\Omega_e^h}
        +\Bigg\langle \frac{ \partial \bfNE}{\partial x_2} \,, \derivev{ \bfG_2^{(k)} } \Bigg\rangle_{\Omega_e^h}
        \\
        \derivee{ \bfR_{Ee}^{(k)} }&=  \Bigg\langle \bfNE
        \, , \bfNE \Bigg\rangle_{\Omega_e^h} \,.
    \end{align}
\end{subequations}
Note that we use the definition of a dot product for matrices as $\langle \bfA_1 
        \, ,  \bfA_2  \rangle_{\Omega_e^h}= \int_{\Omega_e^h} \bfA_1 \cdot \bfA_2 \,\da = \int_{\Omega_e^h} \bfA_1^\top \bfA_2 \, \da$.

We present the vector of conservative quantities $\bfU$ (cf. Equation \eqref{Eq:Uterm}) and its derivatives:
\begin{equation}
\begin{array}{cc}
   \bfU  = \left[
    \begin{array}{c}
           h \\
          h\vel_1  \\
          h\vel_2 
    \end{array}
    \right]  \,,
    &
     \derivev{ \bfU } = \left[
    \begin{array}{ccc}
           1 & 0 & 0\\
          \vel_1 & h & 0 \\
          \vel_2 & 0 & h
    \end{array}
    \right] \bfNV \,.
\end{array}
\end{equation}

Also, the shallow water term $\bfF$ (cf. Equation \eqref{Eq:Fterm}) and its derivatives are detailed in the following:
\begin{eqnarray}
\begin{array}{cc}
    \bfF_1 = \left[ \begin{array}{c} 
    h\vel_1 \\
    h\vel_1^2 +1/2 g_c\,\hth^2 \\
    h\vel_1\vel_2 
    \end{array}\right]  
    \,,
    &
    \derivev{\bfF_1} = \left[ \begin{array}{ccc} 
    \vel_1 & h & 0\\
    \vel_1^2 + g_c\,\hth & 2h\vel_1 & 0\\
    \vel_1\vel_2  & h\vel_2 & h\vel_1
    \end{array}\right]   \bfNV
    \,,
\end{array}
\end{eqnarray}

\begin{eqnarray}
 \begin{array}{cc} 
\bfF_2 = \left[ \begin{array}{c} 
    h\vel_2\\
     h\vel_1\vel_2\\
     h\vel_2^2 +1/2 g_c\,\hth^2 
    \end{array}\right]  
    \,,
    &
    \derivev{\bfF_2} = \left[ \begin{array}{ccc} 
    \vel_2 & 0 & h \\
     \vel_1\vel_2 & h\vel_2 & h\vel_1\\
     \vel_2^2 + g_c\,\hth & 0 & 2h\vel_2
    \end{array}\right]   \bfNV
    \,.
\end{array}
\end{eqnarray}

The derivatives for the flux obtained from the Riemann solver \eqref{Eq:HLLSolver} for $i \in {1,2} $ are computed such that:
\begin{eqnarray}
    \derivev{\nbfF_i} = \begin{dcases}
         \derivev{\bfF^{-}} \quad \textrm{if} \quad S_L\ge 0
        \\
        \frac{S_R\derivev{\bfF_i^{-}} -S_L\derivev{\bfF_i^{+}} + S_LS_R \left(\derivev{\bfUe}-\derivev{\bfUi} \right)}{S_R-S_L} \quad \textrm{if} \quad S_L< 0<S_R
        \\
        \derivev{\bfF_i^{+}} \quad \textrm{if} \quad S_R\le 0
    \end{dcases} \,.
\end{eqnarray}

Furthermore, we present the stress term $\bfQ$ (cf. Equation \eqref{Eq:Qterm}) and its derivatives:
\begin{subequations}
    \begin{align}
        \bfQ_1 &= \left[ \begin{array}{c} 
            0 \\
            h  \sigma_{11} /\rho \\
            h  \sigma_{21} /\rho  
            \end{array}\right] 
            \,,
       \\
        \derivev{\bfQ_1} &= \left[ \begin{array}{ccc} 
            0 & 0 & 0\\
              \sigma_{11} /\rho  & 0 & 0  \\
              \sigma_{21} /\rho  & 0 & 0  
            \end{array}\right]  \bfNV
            \,,
      \\
        \derivee{\bfQ_1} &= \frac{ h}{\rho }\left[ \begin{array}{cccc} 
            0 & 0 & 0& 0 \\
             \derive{\sigma_{11}}{\EE_{11}} & \derive{\sigma_{11}}{\EE_{21}} & \derive{\sigma_{11}}{\EE_{12}} & \derive{\sigma_{11}}{\EE_{22}} \\
             \derive{\sigma_{21}}{\EE_{11}} & \derive{\sigma_{21}}{\EE_{21}} & \derive{\sigma_{21}}{\EE_{12}}  &
             \derive{\sigma_{21}}{\EE_{22}} 
            \end{array}\right] \bfNE
            \,,
    \end{align}
\end{subequations}

\begin{subequations}
    \begin{align}
        \bfQ_2 &= \left[ \begin{array}{c} 
            0 \\
            h  \sigma_{12} /\rho \\
            h  \sigma_{22} /\rho  
            \end{array}\right] 
            \,,
       \\
        \derivev{\bfQ_2} &= 
            \left[ \begin{array}{ccc} 
            0 & 0 & 0\\
              \sigma_{12} /\rho  & 0 & 0  \\
              \sigma_{22} /\rho  & 0 & 0  
            \end{array}\right]  \bfNV
            \,,
      \\
        \derivee{\bfQ_2} &= 
            \frac{ h}{\rho } 
            \left[ \begin{array}{cccc} 
            0 & 0 & 0 & 0\\
            \derive{\sigma_{12}}{\EE_{11}} & \derive{\sigma_{12}}{\EE_{21}} & \derive{\sigma_{12}}{\EE_{12}} & \derive{\sigma_{12}}{\EE_{22}} \\
            \derive{\sigma_{22}}{\EE_{11}} & \derive{\sigma_{22}}{\EE_{21}} & \derive{\sigma_{22}}{\EE_{12}}  &
            \derive{\sigma_{22}}{\EE_{22}}
            \end{array}\right] \bfNE
            \,.
    \end{align}
\end{subequations}

\begin{subequations}
    \begin{align}
    \derivev{\nbfQ} \cdot \bfne &= \frac{1}{2} \left( 
    \derivev{\bfQ_1^+}  \bfne \cdot \bfe_1 +\derivev{\bfQ_2^+}  \bfne \cdot \bfe_2  +\derivev{\bfQ_1^-}  \bfne \cdot \bfe_1 
    + \derivev{\bfQ_2^-}  \bfne \cdot \bfe_2\right) \\
     \derivee{\nbfQ} \cdot \bfne &= \frac{1}{2} \left( 
    \derivee{\bfQ_1^+}  \bfne \cdot \bfe_1 +\derivee{\bfQ_2^+}  \bfne \cdot \bfe_2  +\derivee{\bfQ_1^-}  \bfne \cdot \bfe_1 
    + \derivee{\bfQ_2^-}  \bfne \cdot \bfe_2\right) \,.
\end{align}
\end{subequations}

To obtain the derivatives of the stress for our one-dimensional problem, we derive Equations \eqref{Eq:binghamR}, \eqref{Eq:sigma11D}, and \eqref{Eq:bingham1D} to obtain the following expressions:
\begin{equation}
    \derivees{\bfsigma} = 
     \derivees{\bfsigmaOne}    +     \derivees{\bfsigmaTwo}\,,
\end{equation}
\begin{eqnarray}
    \derivees{\bfsigmaOne}=4\eta \,,
\end{eqnarray}
\begin{eqnarray}
    \derivees{\bfsigmaTwo} = \begin{cases}
       0 & \text{if} \quad \left| \EE \right| \ge \sigmao / \gamma \\
        2 \gamma    & \text{if} \quad \left| \EE \right|< \sigmao / \gamma\,.
    \end{cases}
\end{eqnarray}
For the smooth regularizations of the Bingham stress, we derive Equations \eqref{Eq:sigmaTwoA}, \eqref{Eq:sigmaTwoAA}, \eqref{Eq:sigmaTwoB}, and \eqref{Eq:sigmaTwoC} which are:
\begin{eqnarray}
    \derivees{\bfsigmaTwoa}  &=& 2\sigmao \, \frac{ \gamma  }{\max_{\beta} \left( \gamma\left|\EE \right| - \sigmao\,, 0\right) + \sigmao}
    -
    2\sigmao \, \frac{ \gamma^2 \, \EE\, \signf{\EE} }{ \left(\max_{\beta} \left( \gamma\left|\EE \right| - \sigmao\,, 0\right) + \sigmao\right)^2 } \max_{\beta}' \left( \gamma\left|\EE \right| - \sigmao \,, 0\right)\,,
\end{eqnarray}
\begin{eqnarray}
    \max_{\beta}' \left( x \,, 0\right) =
    \begin{cases}
        1 &\quad \text{if} \quad x\ge \frac{1}{2\beta} \\
        \beta \left(x  + \frac{1}{2\beta} \right) &\quad \text{if} \quad |x|\le \frac{1}{2\beta} \\
        0 &\quad \text{if} \quad x \le -\frac{1}{2\beta}.
    \end{cases}
\end{eqnarray}
\begin{eqnarray}
    \derivees{\bfsigmaTwob} =
    \begin{cases}
        0
        &\quad \text{if} \quad 
        \gamma|\EE|\ge  \sigmao + \frac{1}{2\beta} 
        \\
          2\beta \gamma \left(\sigmao-\gamma |\EE|+\frac{1}{2\beta} \right)  
        &\quad \text{if} \quad \sigmao - \frac{1}{2\beta} \le \gamma|\EE|\le  \sigmao + \frac{1}{2\beta} 
        \\
        2\gamma\, \EE 
        &\quad \text{if} \quad \gamma|\EE|\le  \sigmao - \frac{1}{2\beta}.
    \end{cases}
\end{eqnarray}
\begin{eqnarray}
    \derivees{\bfsigmaTwoc} =
    2 \sigmao \gamma \, \text{sech}^2 \left( \gamma \,\EE \right) \,.
\end{eqnarray}

Additionally, we present the expression for source term $\bfS$ (cf. Equation \eqref{Eq:Sterm}) and its derivatives:
\begin{eqnarray}
    \begin{array}{cc} 
        \bfS = \left[ \begin{array}{c} 
            0 \\ 
            -g_s\, \hth - g_c\,\hth\,  \frac{\partial \htH }{\partial \htx_1} \\ 
         - g_c\,\hth\, \frac{\partial \htH }{\partial \htx_2}
         \end{array}\right] \,,
         &
         \derivev{\bfS} = \left[ \begin{array}{ccc} 
            0 &0 &0\\ 
            -g_s - g_c\,  \frac{\partial \htH }{\partial \htx_1} &0&0\\ 
         - g_c\, \frac{\partial \htH }{\partial \htx_2} &0&0
         \end{array}\right] \bfNV\,.
     \end{array}
\end{eqnarray}

Finally, the term $\bfG$ (cf. Equation \eqref{Eq:Gterm}) is detailed with its derivatives:
\begin{eqnarray}
 \begin{array}{cc} 
    \bfG_1 = \left[ \begin{array}{c} 
            h\vel_1^2 \\
            h\vel_2^2\\
            0 \\
            0 
        \end{array}\right] \,,
    &
    \derivev{\bfG_1} = \left[ \begin{array}{ccc} 
            \vel_1^2 & 2h\vel_1 & 0\\
            \vel_2^2 & 0 & 2h\vel_2\\
            0 & 0 & 0 \\
            0 & 0 & 0 
        \end{array}\right]  \bfNV \,,
\end{array}
\end{eqnarray}
\begin{eqnarray}
 \begin{array}{ccc} 
    \bfG_2 = 
        \left[ \begin{array}{c} 
            0 \\
            0 \\
            h\vel_1^2 \\
            h\vel_2^2
        \end{array}\right] \,,
    &
    \derivev{\bfG_2} = \left[ \begin{array}{ccc} 
            0 & 0 & 0\\
            0 & 0 & 0\\
            \vel_1^2 & 2h\vel_1 & 0\\
            \vel_2^2 & 0 & 2h\vel_2\\
        \end{array}\right]  \bfNV \,,
\end{array}
\end{eqnarray}

\begin{subequations}
    \begin{align}
    \derivev{\nbfG} \cdot \bfne &= \frac{1}{2} \left( 
    \derivev{\bfG_1^+}  \bfne \cdot \bfe_1 +\derivev{\bfG_2^+}  \bfne \cdot \bfe_2  +\derivev{\bfG_1^-}  \bfne \cdot \bfe_1 
    + \derivev{\bfG_2^-}  \bfne \cdot \bfe_2\right)  \,.
\end{align}
\end{subequations}

\bibliographystyle{plain}
\bibliography{SWE}
\end{document}